\documentclass[fleqn,12pt]{article}
\usepackage{amssymb,amsmath,amsfonts}
\usepackage{dsfont}
\usepackage{color}
\usepackage{graphicx}
\usepackage{graphics}
\usepackage{latexsym}
\usepackage{amsmath}
\usepackage{amssymb}
\usepackage{graphics}
\usepackage[dvips]{epsfig}
\usepackage{mathrsfs}

\usepackage[margin=1in]{geometry}

\numberwithin{equation}{section}
\newtheorem{thm}{Theorem}[section]

\newtheorem{lem}[thm]{Lemma}

\newtheorem{rmk}[thm]{Remark}
\newtheorem{defi}[thm]{Definition}

\newcommand{\ea}{\epsilon}
\newcommand{\ta}{\theta}

\newcommand{\da}{\delta}
\newcommand{\la}{\lambda}
\renewcommand{\aa}{\alpha}

\newcommand{\pl}{\partial}
\newcommand{\sa}{\sigma}
\newcommand{\oa}{\omega}
\newcommand{\ga}{\gamma}

\newcommand{\iy}{\infty}

\newcommand{\lt}{\left}
\newcommand{\rt}{\right}
\newcommand{\be}{\begin{equation}}
\newcommand{\bs}{\begin{split}}
\newcommand{\es}{\end{split}}
\newcommand{\ee}{\end{equation}}
\newcommand{\bee}{\begin{equation*}}
\newcommand{\eee}{\end{equation*}}

\newcommand{\ef}{\eqref}

\begin{document}
\begin{center}
\large{ \bf Nonlinear
  Asymptotic Stability of the  Lane-Emden Solutions for the Viscous
  Gaseous Star Problem with Degenerate Density Dependent Viscosities}
\end{center}
\begin{center}
{\small Tao Luo, Zhouping Xin, Huihui Zeng}
\end{center}
\begin{abstract} The nonlinear asymptotic stability of  Lane-Emden solutions is proved in this paper for
spherically symmetric motions of viscous gaseous stars with the density dependent shear and bulk viscosities which vanish at the vacuum, when the adiabatic exponent $\gamma$ lies in the stability
regime $(4/3, 2)$, by establishing the global-in-time regularity uniformly up to the vacuum boundary for the vacuum free boundary problem of the compressible Navier-Stokes-Poisson systems with spherical symmetry,  which ensures the global existence of  strong solutions  capturing the precise
physical behavior that the sound speed is $C^{{1}/{2}}$-H$\ddot{\rm o}$lder continuous  across the vacuum boundary,   the large time asymptotic uniform convergence of the evolving vacuum boundary, density and velocity to those of  Lane-Emden solutions with detailed convergence rates, and the detailed large time behavior of solutions near the vacuum boundary.
The results obtained in this paper extend those in \cite{LXZ} of the authors for the constant viscosities to the case of density dependent viscosities which are degenerate at vacuum states.
\end{abstract}

\section{Introduction}
The evolution of a viscous gaseous star in three spatial dimensions with the boundary being the interface between the gas and vacuum can be modeled by the following free boundary problem
\be\label{0.1}\lt\{ \begin{split}
&   \rho_t  + {\rm div}(\rho {\bf u}) = 0 &  {\rm in}& \ \ \Omega(t), \\
 &    (\rho {\bf u})_t  + {\rm div}(\rho {\bf u}\otimes {\bf u})+{\rm div}\mathfrak{S }= - \rho \nabla_{\bf x} \Psi  & {\rm in}& \ \ \Omega(t),\\
 &\rho>0 &{\rm in }  & \ \ \Omega(t),\\
 &\rho=0 \ \ {\rm and} \ \  \mathfrak{S}{\bf n}={\bf 0}   &    {\rm on}& \  \ \Gamma(t):=\pl \Omega(t),\\
 &    \mathcal{V}(\Gamma(t))={\bf u}\cdot {\bf n}, & &\\
&(\rho,{\bf u})=(\rho_0, {\bf u}_0) & {\rm on} & \ \  \Omega:= \Omega(0).
 \end{split}\rt. \ee
 Here $({\bf x},t)\in \mathbb{R}^3\times [0,\iy)$,  $\rho $, ${\bf u} $, $\mathfrak{S}$ and $\Psi$ denote, respectively, the space and time variable, density, velocity, stress tensor and gravitational potential; $\Omega(t)\subset \mathbb{R}^3$, $\Gamma(t)$, $\mathcal{V}(\Gamma(t))$ and ${\bf n}$ represent, respectively, the changing volume occupied by a fluid at time $t$, moving interface of fluids and vacuum states, normal velocity of $\Gamma(t)$ and exterior unit normal vector to $\Gamma(t)$. The gravitational potential is given by
 \bee\label{potential}\Psi({\bf x}, t)=-G\int_{\Omega(t)} \frac{\rho({\bf y}, t)}{|{\bf x}-{\bf y}|}d{\bf y}, \ \ {\rm satisfying} \  \  \Delta \Psi=4\pi G \rho \  \ {\rm in} \ \  \Omega(t)
\eee
with  the gravitational constant $G$ taken to be unity for convenience.
The  stress tensor takes the form:
\bee\label{0.2}
\mathfrak{S}=pI_3-\lambda_1 \left(\nabla {\bf u}+\nabla {\bf u}^t-\frac{2}{3}({\rm div} {\bf u}) I_3\right)-\lambda_2({\rm div}  {\bf u})I_3, \eee
 where $I_3$ is the $3\times 3$ identical matrix, $p$ is the pressure of the gas, $\lambda_1$ is the shear viscosity, $\lambda_2$ is the bulk viscosity, and $\nabla {\bf u}^t$ denotes the transpose of $\nabla {\bf u}$.
 We consider the polytropic gases for which the equation of state is given by
 \bee\label{1.3}
 p=p(\rho)=K\rho^{\gamma}, \eee
 where $K>0$ is a constant set to be unity for convenience, $\gamma>1$ is the adiabatic exponent.

For a non-rotating gaseous star,   the stable  equilibrium configurations, which minimize the energy among all possible configurations (cf. \cite{liebyau}), are spherically symmetric, called Lane-Emden solutions. Therefore,  spherically symmetric
motions for dynamical problems are important to consider. When the viscosities are positive constants,   the global-in-time spherically symmetric solution  to the free boundary problem \eqref{0.1} and its nonlinear asymptotic stability toward the Lane-Emden solution  were proved in \cite{LXZ} for $4/3<\ga<2$ (the stable index), by establishing the global-in-time regularity uniformly up to the vacuum boundary of  solutions  capturing an interesting behavior called  the {\em physical vacuum}  which states that the sound speed $c=\sqrt{p'(\rho)}$ is $C^{ {1}/{2}}$-H$\ddot{\rm o}$lder continuous near the vacuum boundary (cf. \cite{10,10',jangmas,jm,tpliudamping,13,ya}),
as long as the initial datum is a suitably small perturbation  of the Lane-Emden solution  with the same total mass. The large time asymptotic convergence of the global strong solution, in particular, the convergence of  the vacuum boundary and  the uniform convergence of the density,  to those of the Lane-Emden solution  with  detailed convergence rates as the time goes to infinity are given in \cite{LXZ} when the viscosities are constant. The aim of this work is to extend those results to
the more physically reasonable case that the viscosities depend on the density and vanish at vacuum states. The degeneracy of the viscosities near vacuum is one of the main difficulties which makes such a generalization highly nontrivial.  We assume that the viscosities $\lambda_1$ and $\lambda_2$ are  solely functions of the density in this paper. For simplicity,
we set
 \be\label{start}\lambda_1=\nu_1 \rho^\ta \ \ {\rm and } \ \ \lambda_2=\nu_2 \rho^\ta, \ \ {\rm where} \ \  \nu_1,  \ \   \nu_2 \ \ {\rm and} \ \ \ta \ \ {\rm are} \ \ {\rm positive} \ \  {\rm constants}.\ee
One may check easily that the theorems proved in this article apply to the case that $\lambda_1$ and $\lambda_2$ are positive for $\rho>0$ and have the behavior of $\rho^{\theta}$ as $\rho\to 0$.

In the spherically symmetric setting, that is, $\Omega(t)$ is a ball with the  changing radius $R(t)$,
\begin{equation*}\label{3.1}
 \rho({\bf x}, t) = \rho(r, t) \  \ {\rm and} \ \  {\bf u}({\bf x}, t) =  u(r, t) {\bf x} /r  \ \ {\rm with} \ \   r=|{\bf x}| \in \lt(0, R(t)\rt); \end{equation*}
system  \eqref{0.1} can then be reduced to
\begin{equation}\label{103}
\lt\{\begin{split}
&  (r^2\rho)_t+ (r^2\rho u)_r=0  & {\rm in } & \  \  \lt(0, \  R(t)\rt) , \\
&\rho( u_t +u u_r)+  p_r+ {4\pi\rho}r^{-2}\int_0^r\rho(s,t) s^2ds\\
& \ \  \ \  =\lt[\lt(\frac{4}{3}\la_1+\la_2\rt)\frac{(r^2 u)_r}{r^2} \rt]_r
-4(\la_1)_r\frac{u}{r} & {\rm in } & \  \ \lt(0, \  R(t)\rt), \\
& \rho>0  & {\rm in } & \  \  \lt[0, \  R(t)\rt),\\
& \rho=0   & {\rm for} &  \   \ r=R(t), \\
& \dot R(t)=u(R(t), t) \ \ {\rm with} \ \ R(0)=R_0,  \ \    u(0,t)=0, &    &\\
&  (\rho, u) = (\rho_0, u_0)   & {\rm on } & \ \ (0, \ R_0).
\end{split}\rt.
\end{equation}
The initial domain is taken to be a  ball  $\{0\le r\le R_0\}$. And the initial density is supposed to satisfy
the following  condition:
\be\label{156}
\rho_0(r)>0 \ \ {\rm for} \ \  0\le r<R_0 , \ \ \rho_0(R_0)=0 \ \ {\rm and}  \ \
  -\iy<   \lt(\rho_0^{\ga-1}\rt)_r <0 \  \ {\rm at} \  \ r=R_0.
   \ee
So,
\be\label{physicalvacuum}\rho_0^{\ga-1}(r) \sim R_0-r \ {\rm as~} r {\rm~ close~ to~} R_0, \ee
that is, the initial sound speed is $C^{1/2}$-H${\rm \ddot{o}}$lder continuous across the vacuum boundary, which is called the {\em physical vacuum} for the compressible inviscid flows
(cf. \cite{10,10',jangmas,13,ya}).

The behavior \eqref{156} near the vacuum boundary captures a very interesting feature of the stationary solution of \eqref{103}, $(\rho,u)=(\bar \rho, 0)$, the Lane-Emden solution  (cf. \cite{ch,linss}), corresponding to a non-rotating gaseous sphere in hydrostatic equilibrium. Here $\bar\rho$ solves
\be\label{le} \pl_r(\bar \rho^{\gamma})+ {4\pi}r^{-2} \bar{\rho} \int_0^r \bar \rho(s) s^2ds=0.
\ee
The solutions to \eqref{le} can be characterized by the values of $\gamma$
(cf. \cite{linss}) for given finite total mass $M > 0$. If $\gamma \in (6/5, 2)$, there exists at least one compactly supported solution, and every solution is compactly supported and unique for $\gamma \in (4/3, 2)$.
If $\gamma= 6/5$, the unique solution admits an explicit expression, and it has infinite support.
On the other hand, for $\gamma\in  (1, 6/5)$, there are no solutions with finite total mass.  For $\gamma>6/5$, let $\bar R$ be the radius of the stationary star giving by the  Lane-Emden solution, then it holds (cf. \cite{linss,makino})
\be\label{pvforle} \bar\rho^{\ga-1}(r) \sim \bar R-r \ {\rm as~} r {\rm~ close~ to~} \bar R.\ee

In both astrophysics and the theory of nonlinear PDEs, the problem of nonlinear asymptotic stability of Lane-Emden solutions  is of fundamental importance. As mentioned in \cite{LXZ},  the key to this is to establish the the global-in-time regularity of higher-order  derivatives  of solutions uniformly up to the vacuum boundary, which has been challenging  due to the high degeneracy of system \eqref{103} caused by the behavior \ef{physicalvacuum} near the vacuum boundary, even for the local-in-time existence theory in both the inviscid and viscous compressible flows.   Indeed, the local-in-time well-posedness of smooth solutions to   free boundary problems with physical vacuum was only established recently for compressible inviscid flows (cf. \cite{10, 10', jangmas, jm}, and \cite{LXZ1} which proved a local-in-time well-posedness theory in a new functional space for the three-dimensional compressible Euler-Poisson equations in  spherically symmetric motions). For the vacuum free boundary problem \eqref{103} of the compressible viscous flows featuring the behavior \eqref{physicalvacuum} near the vacuum boundary,  a local-in-time well-posedness theory of strong solutions was established    in \cite{jangnsp} and \cite{duan1},  for the constant and density dependent viscosities, respectively. In \cite{LXZ}, the authors succeeded in establishing such a global-in-time higher-order regularity for the constant viscosities. Compared with the case of constant viscosities, the degeneracy of viscosities at vacuum states makes this a much more challenging task.
In addition, the term $4(\la_1)_r u/{r}=4\nu_1(\rho^{\ta})_r u/{r}$ in the momentum equation, which does not appear in the constant viscosity case, causes serious difficulties in the analysis for the solutions with the behavior \ef{physicalvacuum}.
Indeed, near the vacuum boundary,  $r=R(t)$, this term involves $\ta(R(t)-r)^{-1+ {\ta}/({\ga-1})} $ which is unbounded for $0<\ta<\ga-1$,
where $\ta$ is the constant in \ef{start}.
This difficulty appears even in the study of local existence of strong solutions. In fact, the local existence result of  strong solutions to \ef{103} in \cite{duan1} only holds for $\ta=1$, which avoids the unboundedness of $4\nu_1(\rho^{\ta})_r u/{r}$
for $0<\ta<\ga-1$ when $\ga<2$. In the present work, special care is taken to deal with this term near the vacuum boundary by constructing suitable weights to capture the behavior of velocity near the vacuum boundary to resolve this difficulty.
In addition to this, we also  refine the arguments in \cite{LXZ} substantially in obtaining the $L^\infty$-bound for the first derivative of velocity, for which the weighted $L^2$-estimates of the second derivative are used in \cite{LXZ}, by establishing a pointwise estimate for the
first derivative of velocity near the vacuum boundary only involving the weighted $L^2$-estimates of the first derivative.

When $3/4<\ga<2$, we extend in this article the nonlinear asymptotic stability results in \cite{LXZ} for the constant viscosities to the case of the density dependent
viscosities (1.2) with $0<\ta\le \ga/2$, by proving the existence of  a unique global-in-time  strong solution to \ef{103} and establishing the global-in-time regularity uniformly up to the vacuum boundary, which ensures  the large time asymptotic uniform convergence of the evolving vacuum boundary, density and velocity to those of  the Lane-Emden solution with detailed convergence rates, and detailed large time behaviors of solutions near the vacuum boundary. In particular, we show that every spherical surface moving with the fluid  converges to the sphere  enclosing the same mass inside the domain of the Lane-Emden solution with a uniform convergence rate. Therefore, the convergence of the vacuum boundary $r=R(t)$ as $t\to \infty$ to that of the Lane-Emden solution $r=\bar R$ holds as a consequence.  This also shows that  the large time asymptotic states for the vacuum free boundary problem \ef{103}  are  determined by the initial mass distribution and the total mass. The results obtained in the present work  are among  few results of  global-in-time {\it strong} solutions to vacuum free boundary problems of compressible fluids capturing the singular behavior of \eqref{physicalvacuum}.

Extensive works have been done on the studies of the Euler-Poisson and  Navier-Stokes-Poisson equations with vacuum, especially in recent years. One may find the study of the stability problem of  gaseous stars  in astrophysics literatures (cf.  \cite{ch,weinberg,lebovitz1}).  The linear stability of Lane-Emden solutions was studied in \cite{linss}.    By assuming the existence of global solutions of the Cauchy problem for  the three-dimensional compressible Euler-Poisson equations which has been a major challenge in the theory of fluid dynamics equations,  a conditional nonlinear Lyapunov type
stability theory of stationary solutions  for $\gamma> 4/3$ was established  in \cite{rein} using a variational approach (the same type of nonlinear stability results for rotating
stars were given  by \cite{luosmoller1,luosmoller2}).
Those nonlinear stability results are in the framework of initial value problems in the entire $\mathbb{R}^3$ and involve only Lyapunov functionals which are essentially equivalent to $L^p$-norms of the difference of solutions, where the vacuum boundary cannot be traced. In the framework of free boundary problems for the Euler-Poisson  and Navier-Stokes-Poisson equations, the nonlinear dynamical instability of Lane-Emden solutions for $\gamma\in (6/5, 4/3)$, was proved in \cite{17'} and \cite{jangtice}, respectively. In inviscid flows,  a nonlinear instability for $\gamma= {6}/{5}$ was proved in \cite{jang65};  an instability was identified  for  $\gamma= {4}/{3}$ in \cite{DLYY} that a small
perturbation can cause part of the mass to go off to infinity.

 It should be noted that the existence of  global weak solutions was proved in \cite{fangzhang1} for the  initial boundary value problem reduced from  the  vacuum free boundary problem \ef{103} after using  the Lagrangian mass coordinates,    under the constraint  that $\ga>4/3$, $\ta\in (0,\ga-1)\cap (0, \ga/2]$, and
 $$\frac{2(8+4\alpha-\alpha^2)}{4-4\alpha+\alpha^2} + \frac{2}{3}- \frac{8\sqrt{5+4\alpha-\alpha^2}}{4-4\alpha+\alpha^2} < \frac{\nu_2}{\nu_1}<\frac{2(8+4\alpha-\alpha^2)}{4-4\alpha+\alpha^2} + \frac{2}{3} + \frac{8\sqrt{5+4\alpha-\alpha^2}}{4-4\alpha+\alpha^2}$$
 with $\alpha\in(-1,1)$ being a constant.    In contrast to the strong stability result shown here, for the global weak solutions obtained in \cite{fangzhang1}, only the uniform convergence of the velocity is proved, due to the lack of regularity near the vacuum boundary, and the uniform convergence of the density, in particular, the convergence of the vacuum boundary which is the most interesting part in the study of asymptotic behavior of the free boundary problem, are missing.
Furthermore, our nonlinear asymptotic stability result holds for $4/3<\ga<2$, $0<\ta\le \ga/2$, $\nu_1>0$ and $\nu_2>0$, without the restrictions on $\ta<\ga-1$ and the ratio of $\nu_2/\nu_1$  as in \cite{fangzhang1}.  It should be noted that $\ga-1<\ga/2$ for $\ga<2$.

We conclude the introduction by reviewing some previous works on viscous flows.
It should be noted that there are also other  prior results on free boundary problems involving vacuum for the compressible Navier-Stokes equations besides the ones aforementioned.  One may refer to   \cite{Okada,Okada3,LXY,fangzhang,JXZ,duan,YYZ,yangzhu,JXZ,zhu}  and references therein for the one-dimensional motions concerning global weak solutions. For the spherically symmetric motions, global existence and stability of weak solutions were obtained in \cite{Okada1,OSM}  for gases  surrounding a solid ball (a hard core),  restricted to cut-off domains excluding a neighborhood of the origin;  a global existence of weak solutions containing the origin was established in \cite{GLX} for which the density does not vanish on the boundary and the  viscosities are  density dependent.  The readers may refer to \cite{94, 95} for the local-in-time well-posedness results and \cite{96} for linearized stability results of stationary solutions for a class of free boundary problems of the compressible Navier-Stokes-Poisson equations away from vacuum states.

\section{Main Results}
First, we recall some properties for Lane-Emden solutions. For $\ga\in(4/3, 2)$, it is known that for any given finite positive total mass, there exists a unique solution to equation \eqref{le} whose support is compact (cf. \cite{linss}). Without abusing notations and for convenience, we use $x$ as the variable in the study of Lane-Emden solutions. That means, for any $M\in(0,\infty)$, there exists a unique function $\bar\rho(x)$ such that
\begin{equation}\label{lex1}
\bar\rho_0:=\bar\rho(0)>0, \ \  \bar\rho(x)>0 \ \  {\rm for} \ \ x\in \lt(0,\ \bar R\rt),   \ \  \bar\rho\lt(\bar R\rt)=0, \ \  M=\int_0^{\bar R} 4\pi \bar \rho(s) s^2ds  ;
 \end{equation}
 \begin{equation}\label{newle}
-\iy<\bar\rho_x<0 \ \ {\rm for} \ \  x\in (0,\ \bar R) \ \ {\rm and} \ \ \bar\rho(x) \le \bar\rho_0 \ \  {\rm for} \ \ x\in \lt(0,\ \bar R\rt);
\end{equation}
\begin{equation}\label{rhox}
\left(\bar{\rho}^\ga \right)_x=-x \phi \bar{\rho}, \ \ {\rm where} \  \ \phi:= x^{-3}\int_0^x 4\pi \bar\rho(s)  s^2 ds \in \left[M/{\bar R}^3, \   4\pi \bar\rho_0/3\right] ;
\end{equation}
for a certain finite positive constant $\bar R$ (indeed, $\bar R$ is determined by $M$ and $\ga$). Note that
\begin{equation}\label{ga1x}
\left(\bar{\rho}^{\ga-1} \right)_x= \frac{\ga-1}{\ga}  {\bar \rho}^{-1} \left(\bar{\rho}^{\ga} \right)_x
= -\frac{\ga-1}{\ga}x\phi.
\end{equation}
It then follows from \eqref{lex1} and \eqref{rhox} that $\bar\rho$ satisfies the physical vacuum condition, i.e.,
\begin{equation*}\label{newphy}
 \bar{\rho}^{\ga-1}(x)  \sim  \bar R- x \ \ {\rm as} \ \ x \ \ {\rm~ close~ to~} \bar R.
\end{equation*}
Indeed, there exists a constant $C$ depending on $M$ and $\ga$ such that
\begin{equation}\label{phy}
C^{-1} \lt( \bar R- x \rt) \le \bar{\rho}^{\ga-1}(x)  \le C \lt( \bar R- x \rt), \ \ x\in \lt(0,\ \bar R\rt).
\end{equation}

We adopt a particle trajectory Lagrangian formulation for  \eqref{103} as follows.   Let  $x$ be the reference variable and define the Lagrangian variable $r(x, t)$ by
$$
  r_t(x, t)= u(r(x, t), t) \ \ {\rm for} \  \ t>0  \ \  {\rm and} \ \ r(x,0)=r_0(x), \ \  x\in I:=\lt(0, \bar R\rt)  .
$$
 Here  $r_0(x)$ is the initial position which maps  $ \bar I \to \lt [0, R_0\rt]$  satisfying
\begin{equation}\label{rox} \int_0^{r_0(x)} \rho_0(s) s^2ds = \int_0^x   \bar\rho(s)  s^2 ds, \ \ x\in \bar I, \end{equation}
so that
\begin{equation}\label{choice}
\rho_0(r_0(x))r_0^2(x)r_0'(x)=\bar\rho(x)x^2, \ x\in \bar I.\end{equation}
(Indeed, \ef{rox} means that the initial mass in the ball with the radius $r_0(x)$ is the same as that of the Lane-Emden solution in the ball
with the radius $x$. Then smoothness of $r_0(x)$ at $x=\bar R$ is equivalent to that the initial density $\rho_0$ has the same behavior
near $R_0$ as that of $\bar\rho$ near $\bar R$.)
The choice of  $r_0$  can be described by
\be\label{r000} r_0(x)=\psi^{-1}(\xi(x)), \ \ 0\le x\le \bar R; \ee
where $\xi$ and $\psi$ are one-to-one mappings, defined by
$$\xi:  (0, \bar R) \to (0, M): \ x \mapsto \int_0^x s^2\bar \rho (s)d s  \ \ {\rm and} \ \
  \psi:  (0, R_0) \to (0, M): \ z \mapsto \int_0^z s^2 \rho_0 (s)ds. $$
Moreover $r_0(x)$ is an increasing  function and
the initial total mass has to be the same as that for $\bar \rho$, that is,
\begin{equation}\label{samemass}
\int_0^{R_0}4\pi  \rho_0(s) s^2 ds=\int_0^{r_0(\bar R)}4\pi  \rho_0(s) s^2 ds =
\int_0^{\bar R} 4\pi  \bar \rho (s) s^2 ds=M,
\end{equation}
to ensure that $r_0$ is a diffeomorphism from $\bar I$ to $\lt [0, R_0\rt]$.
It follows from $\eqref{103}_1$, that
\begin{equation}\label{masswithin}
\int_0^{r(x, t)}\rho(s, t)s^2ds=\int_0^{r_0(x)} \rho_0(s) s^2ds, \ \  x\in I.
\end{equation}

Set
$$f(x,t)=\rho(r(x,t),t) \ \ {\rm and} \ \ v(x,t)=u(r(x,t),t).$$
Then  the Lagrangian version of system $\eqref{103}_{1,2}$ can be written on the reference domain $I$ as
\begin{equation}\label{new419}\lt\{ \begin{split}
& (r^2f)_t +r^2f\frac{v_x}{r_x}=0, \\
& f  v_t+ \frac{ (f^{\gamma})_x}{  r_x}+ {4\pi f} {r^{-2}} \int_0^{r_0(x)} \rho_0(s) s^2ds = \frac{1}{r_x } \left[\lt(\frac{4}{3}\la_1+\la_2\rt)\frac{ (r^2v)_x}{r^2  r_x}\right]_x
 -4(\la_1)_x\frac{v}{ r_x r}.
\end{split} \rt.
\end{equation}
Solving $\eqref{new419}_1$ gives   that
$$f(x, t)r^2(x, t)  r_x(x, t)= \rho_0(r_0(x)) r_0^2(x)  r_{0x}(x), \ \  x\in I.$$
Therefore,
$$
f(x, t)=  \frac{x^2\bar \rho(x)}{r^2(x, t)  r_x(x, t)}  \ \ {\rm for} \   \ x\in  I,
$$
 due to   \ef{choice}.
By using \ef{rhox}, the free boundary problem \ef{103} is then reduced to the following initial boundary value problem on a fixed interval $\bar I$:
\begin{equation}\label{419}\lt\{\begin{split}
 & \bar\rho\left( \frac{x}{r}\right)^2  v_t   +  \left[    \left(\frac{x^2}{r^2}\frac{\bar\rho}{ r_x}\right)^\ga    \right]_x +  \frac{x^2}{r^4}    \bar\rho \int_0^x 4\pi \bar\rho  y^2 dy= \mathcal{V}(x,t)
  & {\rm in}  & \  \   I\times (0, T], \\
&  v(0, t)=0    &  {\rm on}& \  \ (0,T],\\
& (r,\ v)(x, 0) = \lt(r_0(x), \  u_0(r_0(x)) \rt)    & {\rm on}& \  \   I \times \{t=0\},
\end{split} \rt.
\end{equation}
where
\bee
\mathcal{V}= \nu \left[\left(\frac{x^2}{r^2}\frac{\bar\rho}{ r_x}\right)^\ta\frac{ (r^2v)_x}{r^2  r_x}\right]_x  -4\nu_1\left[    \left(\frac{x^2}{r^2}\frac{\bar\rho}{ r_x}\right)^\ta\rt]_x\frac{v}{ r}   \  \ {\rm with} \ \    \nu=\frac{4}{3}\nu_1 +\nu_2>0.
\eee
It should be noticed that $\mathcal{V}$ can be rewritten as
\be\label{abc}
\mathcal{V} =-\frac{\nu}{\ta}\lt(\frac{r}{x}\rt)^{-\frac{4\nu_1}{\nu}\ta}
 \lt\{\lt(\frac{r}{x}\rt)^{\frac{4\nu_1}{\nu}\ta} \lt[\bar\rho^\ta \lt(\frac{x^2}{r^2r_x}\rt)^\ta \rt]_x\rt\}_t.
 \ee

A strong solution to  problem \ef{419} is defined as follows.

\begin{defi} $v\in L^{\infty}\lt([0, T]; H^2_{loc}([0, \bar R))\rt)\cap L^{\infty}\lt([0, T];  W^{1, \infty}(I)\rt)$ with
\be\label{r}
r(x, t)=r_0(x)+\int_0^t v(x, s)ds  \ \ {\rm for } \ \  (x, t)\in I\times [0, T]
\ee
 satisfying the initial condition $\ef{419}_3$ is called a strong solution of problem \ef{419} in $[0, T] $, if

1)  $r_x(x, t)>0$ for $(x, t)\in I\times [0, T]$;

2)  $\bar\rho^{ {1}/{2}} v \in C^1([0, T]; L^2(I))$;

3)  $r\in L^{\infty}\lt([0, T]; H^2_{loc}([0, \bar R))\rt)$ and  $ \lt( \bar\rho^{\ga-1/2} (r/x)_x, \  \bar\rho^{\ga-1/2}r_{xx}, \   \bar\rho^{-1/2}\mathcal{V} \rt)\in L^{\infty}([0, T]; L^2(I)) $;

4) $v(0, t)=0$  holds in the sense of $W^{1, \infty}$-trace   for $t\in [0, T]$;

5) $\ef{419}_1$ holds for $(x, t)\in I\times [0, T]$, a.e..
\end{defi}

We are ready to state the main theorem of this paper. Denote
\be\label{lownorm}\begin{split}
\mathscr{E}(t)=   \|(r_x-1, \ v_x)(\cdot,t)\|^2_{L^\iy}
+
\lt\|\bar\rho^{\ga-1/2} \lt( (r/x)_x, \  r_{xx}\rt) (\cdot, t) \rt\|_{L^2}^2
+
\lt\|\bar\rho^{1/2} v_t (\cdot, t) \rt\|_{L^2}^2.
\end{split} \ee

\begin{thm}\label{hhmainthm1} Let $\ga\in(4/3,\ 2)$, $0<\ta\le \ga/2$,  and $\bar\rho$ be the Lane-Emden solution satisfying \ef{lex1}-\ef{rhox}. Assume that the compatibility condition $v(0, 0)=0$ holds and the initial density $\rho_0$ satisfies \eqref{156} and \ef{samemass}.
There exists a  constant  $\bar\da >0$ such that if
$$\mathscr{E}(0)\le \bar\da, $$
then   the  problem \eqref{419}  admits a unique strong solution  in $I\times[0, \iy)$ with
$$\mathscr{E}(t)\le C \mathscr{E}(0) , \ \ \ \  t\ge 0, $$
for  some constant $C$ independent of $t$.
\end{thm}
It should be noted that  $\lt\|\bar\rho^{1/2} v_t (\cdot, 0) \rt\|_{L^2}$ is given in terms of the initial data $(r_0, u_0)$ by   equation $\ef{419}_1$.

For any $t\ge 0$, since $r_x(x, t)>0$ for $x\in \bar I$, $r(x, t)$ defines a diffeomorphism from the reference domain  $\bar I$ to the changing domain $\{0\le r\le R(t)\}$ with the boundary
\be\label{vacuumboundary}
R(t)=r\lt(\bar R ,  t\rt).
\ee
 It also induces a diffeomorphism from the initial domain, $\bar B_{R_0}(0)$,
 to the evolving domain, $\bar B_{R(t)}(0)$,
  for all $t\ge 0$:
$${\bf x}\ne {\bf 0} \in \bar B_{R_0}(0)\to r\lt(r_0^{-1}(|{\bf x}|),  t\rt)\frac{{\bf x}}{|{\bf x}|} \in  \bar B_{R(t)}(0), $$
where $r_0^{-1}$ is the inverse map of $r_0$ defined in \ef{r000}. Here
$$\bar B_{R_0}(0):= \{{\bf x}\in \mathbb{R}^3: |{\bf x}|\le R_0\} \ \  {\rm and} \  \ \bar B_{R(t)}(0):=\{{\bf x}\in \mathbb{R}^3: |{\bf x}|\le R(t)\}. $$
 Denote  the inverse of the map $r(x, t)$ by $\mathcal{R}_t$ for  $t\ge 0$ so that
$$
{\rm if~} \ \  r=r(x, t)\ \ {\rm  for ~} \ \  0\le r\le R(t), \ \   {\rm then~} \ \ x=\mathcal{R}_t(r).$$
For the strong solution $(r, v)$ obtained in Theorem \ref{hhmainthm1}, we set for $0\le r\le R(t)$ and $t\ge 0$,
\be\label{solution}
  \rho(r, t)=\frac{x^2\bar \rho(x)}{r^2(x, t)  r_x(x, t)}  \ \  {\rm and} \ \
   u(r, t)=v(x, t) \ \  {\rm with} \ \     x=\mathcal{R}_t(r).  \ \
   \ee
Then the triple $(\rho(r,t), u(r,t), R(t))$  ($t\ge 0$) defines a global strong solution to the free boundary problem \ef{103}. Furthermore, the strong nonlinear asymptotic stability of the Lane-Emden solution can be stated as follows.

For any $\iota\in \lt(0,\ (2\ga-2-\ta)/8 \rt]$, we set
\be\label{alpha} \alpha=\min\{\ga-1+\ta, \ 2(\ga-1)\}-\iota  , \ee
\be\label{beta} \beta=1+ ({\aa -\iota}) /({\ga-\ta}),\ee
\be\label{varsa}\varsigma=\min\lt\{1, \ \   2^{-1}{\beta } + 2^{-1}({\beta-1})\min\lt\{ 1, \ ({\ga-\ta})/{\aa}\rt\}  \rt\}.\ee
It should be noted that $1<\beta<3$ and $0<\aa-\ta<(\ga-1)$ for  $4/3<\ga<2$  and $0<\ta\le \ga/2$.

\begin{thm}\label{mainthm2} Under the assumptions in Theorem \ref{hhmainthm1}, the triple $(\rho, u, R(t))$   defined by
\ef{vacuumboundary} and \ef{solution} is the unique global strong solution to the free boundary problem \eqref{0.1} satisfying
$R\in W^{1, \infty}( [0, \ +\infty) ).$  Moreover, the solution satisfies the following  estimates:
 for any $\iota\in \lt(0,\ (2\ga-2-\ta)/8 \rt]$, $l\in (0,1)$ and $b\in [0, 2-\gamma]$, there exists  positive constants $C_\iota$ and $C_{\iota,l}$ independent of $x$ and $t$ such that for all $t\ge 0$,
\be\label{decayofr}
 (1+t)^{\frac{\ga-1+\aa-\theta}{\ga+\aa-\ta}\beta}  |r(x, t)-x|^2 \le C_\iota \mathscr{E}(0)  , \ \ x\in I,
\ee
\be\label{decayforrho}
 (1+t)^{\min\lt\{\frac{\beta}{2}-\frac{\max\{0, \ 3\ga-5+2b\}}{2(\ga-\ta)} , \   \beta-1\rt\} } \bar\rho^{-b}(x)  |\rho(r(x,t),t)-\bar\rho(x)|^2\le C_\iota \mathscr{E}(0)  , \ \ x\in I,
\ee
\be\label{decayofrx}
 (1+t)^{\beta-1} \lt( | r_x(x,t) - 1 |^2 + |x^{-1}r(x,t)-1|^2  \rt) \le C_{\iota, l} \mathscr{E}(0)  , \ \ x\in [0,l], \ee
\be\label{decayofv}
(1+t)^{\beta/2} |u(r(x, t), t)|^2 \le C_\iota \mathscr{E}(0)  , \ \ x\in I, \ee
\be\label{deforvx}
(1+t)^{a}  \lt( | u_r(r ,t)  |^2  + |r^{-1}u(r ,t) |^2  \rt) \le C_{\iota } \mathscr{E}(0)  , \ \ x\in I. \ee
Here $a>0$ is given by
\begin{align*}
a=   \min  \lt\{ \beta-1, \ \  \frac{\ga-1+\aa-\theta}{\ga+\aa-\ta}\beta, \ \ \frac{3\beta+\varsigma}{4}- \frac{\beta-\varsigma}{4 \aa} \max\{0, \ 4\ta-4(\ga-1)-\aa\}, \rt.\\
 \lt.\beta -\frac{\beta-\varsigma}{2\aa}\max\{0, \ 2\ta-2(\ga-1) \}, \ \
\frac{\beta}{2}-\frac{\beta}{2(\ga+\aa-\ta)}\max\{0, \ 4\ta-\ga-1\}\rt\} .
\end{align*}
\end{thm}

\begin{rmk} The uniform convergence of $r(x,t)$ to $x$ is given in \ef{decayofr}. This implies that  every spherical surface moving with the fluid  converges to the sphere  enclosing the same mass inside the domain of the Lane-Emden solution with a uniform convergence rate, and  the large time asymptotic states for the vacuum free boundary problem \ef{103}  are  determined by the initial mass distribution and the total mass. In particular, this also implies the convergence of the vacuum boundary:
$$|R(t)-\bar R|\le  C_\iota  (1+t)^{-\frac{\ga-1+\aa-\theta}{2(\ga+\aa-\ta)}\beta}   \sqrt {\mathscr{E}(0)}.$$
Moreover, a better decay rate for $|r-x|$ away from the vacuum boundary can be given by
$$
 (1+t)^{\frac{3\ga-2+2(\aa-\theta)}{2(\ga+\aa-\ta)}\beta-\frac{1}{2}}    |r(x, t)-x|^2 \le C_{\iota,l} \mathscr{E}(0)  , \ \ x\in [0,l], 0<l<1, \ \ t\ge 0.
$$
(See \ef{decayofr-est} for the details.)
The uniform convergence of the density and the velocity to  those of the Lane-Emden solution with  uniform convergence rates are given  by \ef{decayforrho} and \ef{decayofv}, respectively.
The estimate  \ef{decayforrho}  yields not only the uniform convergence for large time but also the behavior of the density  near the vacuum boundary since $\ga<2$.
\end{rmk}

\section{Proof of Theorem \ref{hhmainthm1}}
In this section, we derive some priori estimates under the following  a priori  assumptions. Let $v$ be a strong solution to \ef{419} in the time interval $[0 , T]$ with
$$r(x, t)=r_0(x)+\int_0^t v(x, \tau)d\tau, \ \ (x,t)\in [0, \ 1]\times[0,T],$$
 satisfying  the following {\it a priori}  assumption:
\begin{equation}\label{rx}
   \left|r_x-1\right|+ \left|r/x -1 \right|   \le \ea_0 \ \ {\rm for} \ \
  (x,t)\in I\times [0, T],
\end{equation}
where $\ea_0\in (0, 1/2]$ is a sufficiently small but fixed constant (indeed, $\ea_0$ is required to be less than a  constant depending only on $\ga$); and
\begin{equation}\label{vx}
 \lt|v_x\rt|+ \lt|v/x \rt|  \le 1 \ \ {\rm for} \ \
  (x,t)\in I\times [0, T].
\end{equation}
 In particular, it holds that
\begin{equation}\label{Liy}
   {1}/{2}\le r_x,   {r}/{x} \le  {3}/{2} \ \ {\rm for}
 \ \ (x,t)\in I\times [0, T].
\end{equation}
(Indeed, the a priori  assumption \ef{rx} and \ef{vx} will be verified  in  Lemmas \ref{lem3} and \ref{lem1521}.)

The following  Hardy's  inequalities will be used often to derive the a priori estimates, whose proof can be found in \cite{KM} or
\cite{17'}.
\begin{lem}\label{Hardy} {\rm(Hardy's inequalities)} Let $k>1$  be a given real number and $g$ be a function. If g satisfies that $\int_0^{1/2}x^k(g^2+ g_x^2)dx<\infty$, then it holds that
\be\label{hd1}
\int_0^{1/2}x^{k-2}  g^2 dx\le\int_0^{1/2}x^k (g^2+ g_x^2)dx<\infty.
\ee
Similarly, if $\int_{1/2}^1(1-x)^k(g^2+ g_x^2)dx<\infty$,
then
\be\label{hd2}
\int_{1/2}^1(1-x)^{k-2}  g^2 dx\le\int_{1/2}^{1}(1-x)^k (g^2+ g_x^2)dx<\infty.
\ee
\end{lem}

\subsection{Lower-order estimates}\label{sect3.1}
To derive the lower-order weighted energy estimates, one rewrites equation $\eqref{419}_1$ as
 \begin{equation}\label{nsp}\begin{split}
 \bar\rho\left( \frac{x}{r}\right)^2  v_t   +  \left[    \left(\frac{x^2}{r^2}\frac{\bar\rho}{ r_x}\right)^\ga    \right]_x - \frac{x^4}{r^4}  \left(\bar{\rho}^\ga\right)_x   =  \nu \left[\left(\frac{x^2}{r^2}\frac{\bar\rho}{ r_x}\right)^\ta\frac{ (r^2v)_x}{r^2  r_x}\right]_x  -4\nu_1\left[    \left(\frac{x^2}{r^2}\frac{\bar\rho}{ r_x}\right)^\ta\rt]_x\frac{v}{ r},
\end{split}
\end{equation}
where \eqref{rhox} has been used. Here
$\nu={4}\nu_1/3+\nu_2$.

We outline the analysis for the lower-order estimates here and give some motivations for the proofs. In Lemma \ref{lem2}, we derive the decay estimates for the zeroth-order energy:
$$\int x^2\lt\{ \bar{\rho} v^2 + \bar{\rho}^\ga\left[\left(\frac{r}{x}-1\right)^2  + \left(r_x-1\right)^2 \right] \rt\}(x,t)dx
$$
and the boundedness of the weighted energy:
$$\int x^2 \bar\rho^{\ta}\lt[ \lt(r_x-1\rt)^2 +\lt(\frac{r}{x}-1\rt)^2 \rt] (x,t)dx,$$
by constructing suitable nonlinear weighted functionals.
In Lemma \ref{lem4}, the decay estimates for weighted norms of the higher time-derivatives (than those in Lemma \ref{lem2}) are derived.  Further regularity near the vacuum  boundary and decay estimates are obtained in  Lemmas \ref{lem3.4} and
\ref {lem3.5}, by showing the boundedness of
      $$\int  x^2  \bar\rho^{\theta-\aa} \left[\left(\frac{r}{x}-1\right)^2  + \left(r_x-1\right)^2 \right] (x,t) dx, $$
and the decay of the following weighted quantities:
$$\int  x^2 \bar\rho^{\theta}\left[\left(\frac{r}{x}-1\right)^2  + \left(r_x-1\right)^2 \right] (x,t) dx,$$
$$\int x^2\lt\{ \bar{\rho} \lt( v^2 +  v_t^2 \rt)+ \bar{\rho}^\ga\left[\left(\frac{r}{x}-1\right)^2  + \left(r_x-1\right)^2 \right] \rt\}(x,t)dx,$$
$$\int \lt(v^2 + x^2 \bar\rho^{\theta} v_{x}^2 \rt) (x,t) dx    \ \ {\rm and} \ \   \int\bar\rho^{\theta- \aa/2}(x^2v_{x}^2+v^2) (x,t) dx,$$
where $\alpha$ is given  by \ef{alpha}.  Those two lemmas play a crucial role to the global regularity
uniformly up to the vacuum boundary. In the proof of these two lemmas, we use the multipliers:
$$\int_0^x \bar\rho^{- \alpha}(y)\lt(r^3-y^3\rt)_y dy \ \  {\rm and} \ \ \int_0^x \bar\rho^{-\aa}(y)(r^2 v)_ydy.$$
In Lemma \ref{lem3.5}, we also obtain the  decay estimates of
$$ \int (r(x, t)-x)^2dx \ \ {\rm and} \ \ \int x^2 (r_x(x, t)-1)^2dx .$$

With the above lower-order estimates in hand, we can bound
 $$\sup_{x\in I} \lt( x^3 |r_x(x, t)-1|^2 \rt) \ \ {\rm and} \ \ \sup_{x\in I}  \lt( x^3|v_x(x, t)|^2 \rt), $$
 and estimate the decay of
 $$\sup_{x\in I} \lt(x |r(x, t)-x|^2 \rt) \ \  {\rm and} \ \  \sup_{x\in I} \lt( x |v(x, t)|^2\rt) $$
in Lemma \ref{lem3}, which in particular give the uniform  estimates away from the origin. The key idea is to integrate \eqref{nsp} both in $x$ and $t$ to derive an ODE  for the quantity $Z$ defined in \ef{Z} whose leading term is
$r_x-1$ to obtain the  estimate for $  x^3|r_x(x, t)-1|^2$. In turn, the leading term for $\pl_t Z$ is $v_x$ so that the bound for $\pl_t Z$ yields the bound for $x^3 |v_x(x,t)|^2$. In the proof of this lemma, the balance of the pressure and gravitational filed plays an important role.

It worths pointing out the difficulties in the lower-order estimates for the case of density dependent viscosities. For example, in the proof of Lemma \ref{lem2}, it is shown that the quantity
 \begin{align}\label{6/21}
& (1+t)  \int x^2\lt\{ \bar{\rho} v^2 + \bar{\rho}^\ga\left[\left(\frac{r}{x}-1\right)^2  + \left(r_x-1\right)^2 \right] \rt\}(x,t)dx  \notag\\
 \le & C  \lt( \textrm{initial data} + \int_0^t  \int  x^2\bar{\rho}^\ga \left[\left(\frac{r}{x}-1\right)^2  + \left(r_x-1\right)^2 \right]dxds  \rt) .
\end{align}
In order to obtain the decay of the zeroth-order energy, we show that the double integral on the second line of \ef{6/21} can be bounded by  the initial data and   the term
\begin{equation}\label{xyzt}
C \int_0^t \int\lt[ \bar{\rho}^\ta x\lt(|r_x-1|+\lt|\frac{r}{x}-1\rt|\rt)|v|+  \bar{\rho}^\ta \lt|r-x\rt||xv_x|\rt]dxds,
\end{equation}
using the multiplier $r^3-x^3$ motivated by the virial equations adopted in the study of stellar dynamics and equilibriums  (cf. \cite{lebovitz2, tokusky})  to detect the detailed balance between the pressure and the self-gravitation.
The  term  \ef{xyzt} appears due to the dependence of the viscosities on the density. This is in sharp contrast to the case of constant viscosities  for which the multiplier $r^3-x^3$ matches  the viscosities well and  the double integral
on the second line of \ef{6/21}
is bounded directly by the initial data. It is quite subtle  to bound  \ef{xyzt}
due to the degeneracy of the viscosities at vacuum.
We choose a cut-off function deliberately whose  effective length is a small positive number $\delta$
to localize both near the vacuum boundary and the origin, so that \ef{xyzt}  can be  bounded by
$$C \da^\frac{1}{4}  \int_0^t \int   x^2\bar{\rho}^\ga  \left[\left(\frac{r}{x}-1\right)^2  + \left(r_x-1\right)^2 \right] dx  ds$$
and other terms.
The desired estimates are then obtained by choosing $\delta$ small.
In the estimates of $ x^3 |r_x(x, t)-1|^2$ and $x^3 |v_x(x, t)|^2$ in Lemma \ref{lem3}, we derive an ODE for the quantity $Z$ defined in \ef{Z}, while in the case of constant viscosities, an ODE is also derived for the quantity $ ({\bar\rho(x)})^{-1}{\rho(r(x, t), t)}$ which is simpler than that for $Z$ and  can be solved explicitly in some sense. The ODE for $Z$ is more involved and harder
to solve. We have to identify the leading terms and estimate the error terms.

In the following, we give the details and analysis outlined above.

\begin{lem}\label{lem2} Let $\ta\in (0, 1]$. Suppose that \ef{rx} holds  for suitably small constant $\ea_0$. Then,
\begin{equation}\label{lem2est'}\begin{split}
\mathfrak{D}(t) & + \int_0^t \int x^2\bar{\rho}^\ga\left[\left(\frac{r}{x}-1\right)^2  + \left(r_x-1\right)^2 \right] dxds \\
& +   \int_0^t (1+s) \int  \bar\rho^\ta  \left(x^2 v_{x}^2 + v^2  \right)dx ds \le C \mathfrak{D}(0),  \ \  t\in[0,T],
\end{split}
\end{equation}
where
\begin{equation*}\label{}\begin{split}
  \mathfrak{D}(t)   = &  (1+t)  \int x^2\lt\{ \bar{\rho} v^2 + \bar{\rho}^\ga\left[\left(\frac{r}{x}-1\right)^2  + \left(r_x-1\right)^2 \right] \rt\}(x,t)dx\\
 &+  \int x^2 \bar\rho^{\ta}\lt[ \lt(r_x-1\rt)^2 +\lt(\frac{r}{x}-1\rt)^2 \rt] (x,t)dx  .
\end{split}
\end{equation*}
\end{lem}
{\em Proof}. The proof consists three steps.

{\em Step 1}. In this step, we prove that
for $0<\ta\le 1$ and $\ga>4/3$,
\begin{equation}\label{5-0}\begin{split}
  &(1+t)  \int x^2\lt\{ \bar{\rho} v^2 + \bar{\rho}^\ga\left[\left(\frac{r}{x}-1\right)^2  + \left(r_x-1\right)^2 \right] \rt\}(x,t)dx
\\&+   \int_0^t (1+s) \int \bar{\rho}^\ta   \left(x^2 v_{x}^2 + v^2  \right)  dx ds \\
  \le & C \int \eta(x,0) dx  +C  \int_0^t  \int  x^2\bar{\rho}^\ga \left[\left(\frac{r}{x}-1\right)^2  + \left(r_x-1\right)^2 \right]dxds,
\end{split}
\end{equation}
where
\begin{equation*}\label{etadefn}
{\eta}(x,t)=\frac{1}{2} x^2 \bar{\rho} v^2  + x^2\bar{\rho}^\ga\left[\frac{1}{\ga-1}\left(\frac{x}{r}\right)^{2\ga-2}\left(\frac{1}{r_x}\right)^{\ga-1}
+\left(\frac{x}{r}\right)^{2}r_x - 4 \frac{x}{r}  - \frac{4-3\ga}{\ga-1}\rt].
\end{equation*}

It follows from  \eqref{nsp} and   the boundary condition \eqref{lex1} that, for any $\ell\ge 0$
\begin{align}\label{hz0}
& \frac{d}{dt}\lt\{(1+t)^{\ell}\int  {\eta}(x,t) dx\rt\}  +(1+t)^{\ell} \int \left(\frac{x^2}{r^2}\frac{\bar\rho}{ r_x}\right)^\ta\lt[\nu \frac{ (r^2v)_x}{r^2  r_x}  \left(r^2 v\right)_x  - 4 \nu_1  \left(r v^2 \right)_x \rt]dx\notag\\
&=\ell (1+t)^{\ell-1}\int  {\eta}(x,t) dx.
\end{align}
Each term in the equation above can be estimated as follows. First,  using the Taylor expansion, one may verify that for $\gamma>4/3$,
\begin{equation}\label{etalower}
{\eta}(x,t) \ge \frac{1}{2} x^2 \bar{\rho} v^2 + \frac{3\ga-4}{4} x^2\bar{\rho}^\ga\left[2\left(\frac{r}{x}-1\right)^2  + \left(r_x-1\right)^2 \right],
\end{equation}
\begin{equation}\label{etaup}
{\eta}(x,t) \le \frac{1}{2} x^2 \bar{\rho} v^2 + C(\ga)x^2\bar{\rho}^\ga \left[\left(\frac{r}{x}-1\right)^2  + \left(r_x-1\right)^2 \right],
\end{equation}
provided \ef{rx} holds for a suitably small constant $\epsilon_0$, where $C(\ga)$ is a positive constant depending on $\ga$.
Also,
 $$\nu \frac{ (r^2v)_x}{r^2  r_x}  \left(r^2 v\right)_x  - 4 \nu_1  \left(r v^2 \right)_x
\ge   3 \sa  \left( \frac{r^2}{r_x}v_x^2+ 2r_x v^2  \right),
$$
where $\sa=\min\left\{2\nu_1/3,\ \nu_2 \right\}$.
Integrate \ef{hz0} over $[0, t]$ to obtain, by virtue of \ef{etalower} and \ef{etaup}, that for $0<\ta\le 1$,
\begin{equation}\label{5-0-1}\begin{split}
  &(1+t)^{\ell}  \int x^2\lt\{ \bar{\rho} v^2 + \bar{\rho}^\ga\left[\left(\frac{r}{x}-1\right)^2  + \left(r_x-1\right)^2 \right] \rt\}(x,t)dx
\\&+   \int_0^t (1+s)^{\ell} \int \bar{\rho}^\ta   \left(x^2 v_{x}^2 + v^2  \right)  dx ds \\
  \le & C \int \eta(x,0) dx  +C\ell \int_0^t  (1+s)^{\ell-1} \int x^2 \lt\{  \bar\rho^{\ta} v^2+ \bar{\rho}^\ga \left[\left(\frac{r}{x}-1\right)^2  + \left(r_x-1\right)^2 \right]\rt\}dxds.
\end{split}
\end{equation}
Setting $\ell=0$ in \ef{5-0-1}, leads to
 \begin{equation}\label{lem1est}\begin{split}
  \int x^2\lt\{ \bar{\rho} v^2 + \bar{\rho}^\ga\left[\left(\frac{r}{x}-1\right)^2  + \left(r_x-1\right)^2 \right] \rt\}(x,t)dx &
\\ +   \int_0^t  \int \bar{\rho}^\ta   \left(x^2 v_{x}^2 + v^2  \right)  dx ds
  & \le   C \int \eta(x,0) dx.
\end{split}
\end{equation}
Letting $\ell=1$ in \ef{5-0-1} and using \ef{lem1est} prove \ef{5-0}.

{\em Step 2}. In this step, we prove that for $\ta\in (0,1]$,
\begin{equation}\label{th2}\begin{split}
&\int   x^2 \bar\rho^{\ta}\lt[ \lt(r_x-1\rt)^2 +\lt(\frac{r}{x}-1\rt)^2 \rt] dx
 + \int_0^t \int x^2\bar{\rho}^\ga\left[2\left(\frac{r}{x}-1\right)^2  + \left(r_x-1\right)^2 \right] dxds
\\
& \le   C \mathfrak{D}(0) + C \int_0^t \int\lt[ \bar{\rho}^\ta x\lt(|r_x-1|+\lt|\frac{r}{x}-1\rt|\rt)|v|+  \bar{\rho}^\ta \lt|r-x\rt||xv_x|\rt]dxds.
\end{split}
\end{equation}

It is easy to check, in view of \eqref{nsp} and \eqref{lex1}, that
 \begin{equation}\label{dec}\begin{split}
  \int     \bar{\rho}^\ga \left\{       \left[\frac{x^4}{r^4} \left(r^3-x^3\right)\right]_x  - \left(\frac{x^2}{r^2 r_x} \right)^\ga    \left(r^3-x^3\right)_x   \right\} dx =    -\int x^3 \bar\rho   v_t \left( \frac{r}{x}-\frac{x^2}{r^2}\right)dx\\
- \nu  \int  \bar{\rho}^\ta  \left(\frac{x^2}{r^2 r_x}\right)^\ta \frac{ (r^2v)_x}{r^2  r_x}\left(r^3-x^3\right)_x dx + 4\nu_1 \int \bar{\rho}^\ta  \left(\frac{x^2}{r^2 r_x}\right)^\ta \lt[\frac{v}{ r}(r^3-x^3 )\rt]_x   dx   .
\end{split}
\end{equation}
It yields from simple calculations that
\begin{equation*}\label{}\begin{split}
 & x^{-2}\left\{       \left[\frac{x^4}{r^4} \left(r^3-x^3\right)\right]_x  - \left(\frac{x^2}{r^2 r_x} \right)^\ga    \left(r^3-x^3\right)_x   \right\} \\
  =& 3\left(\frac{x^2}{r^2r_x}\right)^{\ga}   -3\left(\frac{x^2}{r^2r_x}\right)^{ \ga-1 } -   \left(\frac{x}{r}\right)^2r_x
  +4\left(\frac{x}{r}\right)^5 r_x - 7\left(\frac{x}{r}\right)^4 + 4 \frac{x}{r},
\end{split}
\end{equation*}
which, together with Taylor's expansion, gives
\begin{equation}\label{dec1}\begin{split}
& \int     \bar{\rho}^\ga \left\{       \left[\frac{x^4}{r^4} \left(r^3-x^3\right)\right]_x  - \left(\frac{x^2}{r^2 r_x} \right)^\ga    \left(r^3-x^3\right)_x   \right\} dx \\
 \ge & \frac{3(3\ga-4)}{2} \int x^2\bar{\rho}^\ga\left[2\left(\frac{r}{x}-1\right)^2  + \left(r_x-1\right)^2 \right] dx ,
\end{split}
\end{equation}
provided \ef{rx} holds for small $\epsilon_0$.
Simple calculations show that  for $\ta\in(0,1]$,
\begin{equation}\label{dec3}\begin{split}
\nu  \int  \bar{\rho}^\ta  \left(\frac{x^2}{r^2 r_x}\right)^\ta \frac{ (r^2v)_x}{r^2  r_x}\left(r^3-x^3\right)_x dx
=3\nu \frac{d}{dt}\int x^2 \bar{\rho}^\ta  \eta_0(x,t) dx,
 \end{split}
\end{equation}
where
\bee\label{}\begin{split}
\eta_0(x,t)=\lt\{\begin{split}  &\frac{1}{\ta} \left(\frac{x^2}{r^2 r_x}\right)^{\ta}-\frac{1}{\ta-1}\left(\frac{x^2}{r^2 r_x}\right)^{\ta-1}
 + \frac{1}{\ta(\ta-1)}  \ \   & {\rm for }  \ \  & 0<\ta<1, \\
& \frac{x^2}{r^2 r_x}-\ln \lt(\frac{x^2}{r^2 r_x}\rt)-1 \ \  & {\rm for} \ \  & \ta=1.\end{split} \rt.
\end{split} \eee
Moreover,
\begin{equation}\label{dec2}\begin{split}
\int x^3 \bar\rho   v_t \left( \frac{r}{x}-\frac{x^2}{r^2}\right)dx =\frac{d}{dt}\int x^3 \bar\rho   v  \left( \frac{r}{x}-\frac{x^2}{r^2}\right)dx - \int x^2 \bar\rho   v^2  \left( 1+2 \frac{x^3}{r^3}\right)dx\end{split}
\end{equation}
and
 \begin{equation}\label{dec4}\begin{split}
&4\nu_1 \int \bar{\rho}^\ta  \left(\frac{x^2}{r^2 r_x}\right)^\ta \lt[\frac{v}{ r}(r^3-x^3 )\rt]_x   dx  \\
\le & C \int \bar{\rho}^\ta x\lt(|r_x-1|+\lt|\frac{r}{x}-1\rt|\rt)|v| dx +C\int \bar{\rho}^\ta \lt|r-x\rt||xv_x|dx.
\end{split}\ee
It follows from \ef{dec}-\ef{dec4} that
\begin{equation}\label{th1}\begin{split}
& \frac{d}{dt} \int \lt\{ 3\nu x^2 \bar{\rho}^\ta  \eta_0
  + x^3 \bar\rho   v  \left( \frac{r}{x}-\frac{x^2}{r^2}\right) \rt\}dx\\
 &+ \frac{3}{2}(3\ga-4)\int x^2\bar{\rho}^\ga\left[2\left(\frac{r}{x}-1\right)^2  + \left(r_x-1\right)^2 \right] dx
\\
 \le &  C \int \lt\{ \bar{\rho}^\ta x\lt(|r_x-1|+\lt|\frac{r}{x}-1\rt|\rt)|v| + \bar{\rho}^\ta \lt|r-x\rt||xv_x| + x^2 \bar\rho   v^2 \rt\}dx  .
\end{split}
\end{equation}

Using the Taylor expansion, shows that for small $\epsilon_0$ in \ef{rx},
\begin{equation*}\label{}\begin{split}
 \eta_0(x,t)
  \ge & \frac{1}{2}\left(\frac{x^2}{r^2 r_x}-1\right)^2 - C \ea_0\left(\frac{x^2}{r^2 r_x}-1\right)^2   \ge   \frac{1}{4}\left(\frac{x^2}{r^2 r_x}-1\right)^2\\
\ge & \frac{1}{4}\lt[\lt(r_x-1\rt) + 2\lt(\frac{r}{x}-1\rt) \right]^2 - C \ea_0 \lt[\lt(r_x-1\rt)^2 + \lt(\frac{r}{x}-1\rt)^2\rt].
 \end{split}
\end{equation*}
Note that
\begin{equation*}\label{}\begin{split}
&\int x^2 \bar\rho^{\ta}\lt[\lt(r_x-1\rt) + 2\lt(\frac{r}{x}-1\rt) \right]^2dx\\
=&\int x^2 \bar\rho^{\ta}\lt[ \lt(r_x-1\rt)^2 +4\lt(\frac{r}{x}-1\rt)^2 \rt]dx -2\int \lt(x\bar\rho^{\ta}\rt)_x\lt({r}-x\rt)^2dx\\
= &\int x^2 \bar\rho^{\ta}\lt[ \lt(r_x-1\rt)^2 +2\lt(\frac{r}{x}-1\rt)^2 \rt]dx + 2 \frac{\ta}{\ga}\int x^4\phi \bar\rho^{\ta-(\ga-1)}\lt(\frac{r}{x}-1\rt)^2dx
 \end{split}\end{equation*}
where \ef{rhox} has been used. We then have, noting \ef{rhox} again, that
\be\label{cc1}\begin{split}
 C^{-1} \int x^2 \bar\rho^{\ta}\lt[ \lt(r_x-1\rt)^2 +\lt(\frac{r}{x}-1\rt)^2 \rt]dx
&\le \int x^2 \bar{\rho}^\ta  \eta_0dx\\
&\le C \int x^2 \bar\rho^{\ta}\lt[ \lt(r_x-1\rt)^2 +\lt(\frac{r}{x}-1\rt)^2 \rt]dx  .
 \end{split}\ee
Therefore, \ef{th2} follows from   \ef{th1}, the Cauchy inequality and \ef{lem1est}.

{\em Step 3}. In this step, we prove that
\begin{equation}\label{newth3}\begin{split}
  \mathscr{F}(t) \le & C \mathfrak{D}(0) +\int_0^t(1+s)^{- \frac{\ga-\ta}{\ga-\ta- (\ga-1)/2}} \int x^2 \bar{\rho}^{\ta} \lt(|r_x-1|^2+\lt|\frac{r}{x}-1\rt|^2\rt) dxds,
\end{split}
\end{equation}
where
\begin{equation*}\label{}\begin{split}
 \mathscr{F}(t) = & \mathfrak{D}(t)   + \int_0^t \int x^2\bar{\rho}^\ga\left[\left(\frac{r}{x}-1\right)^2  + \left(r_x-1\right)^2 \right] dxds
 \\
 & +   \int_0^t (1+s) \int  \bar\rho^\ta  \left(x^2 v_{x}^2 + v^2  \right)dx ds .
\end{split}
\end{equation*}
If \ef{newth3} is true, then \ef{lem2est'} follows from Grownwall's inequality.

It follows from $\ef{5-0} +k\ef{th2}$ with a suitably large constant $k$ that
\begin{equation}\label{th3}\begin{split}
  \mathscr{F}(t)
  \le  C \mathfrak{D}(0) + C \int_0^t \int\lt[ \bar{\rho}^\ta x\lt(|r_x-1|+\lt|\frac{r}{x}-1\rt|\rt)|v|+  \bar{\rho}^\ta \lt|r-x\rt||xv_x|\rt]dxds.
\end{split}
\end{equation}
Next, we  estimate the last  term on the right-hand side of \ef{th3}. When $A\ge 0$,  it follows from  Young's inequality that
for $0<\oa<({\ga-\ta})/({\ga-1})$,
\begin{align*}
  (1+s)^{-1}\bar\rho^{\ta+\oa(\ga-1)} A &=
\lt\{(1+s)^{-1}\lt(\bar\rho^\ta A\rt)^{\frac{\ga-\ta-\oa(\ga-1)}{\ga-\ta}}\rt\}\lt\{\lt(\bar\rho^\ga A \rt)^{\frac{\oa(\ga-1)}{\ga-\ta}}\rt\}\\
&\le(1+s)^ {- \frac{\ga-\ta}{\ga-\ta-\oa(\ga-1)}} \bar\rho^{\ta}A+\bar{\rho}^{\ga}A.
\end{align*}
Then, for any $0<\oa<({\ga-\ta})/({\ga-1})$,
\begin{equation}\label{hen}\begin{split}
&\int_0^t(1+s)^{-1}\int x^2 \bar{\rho}^{\ta+\oa(\ga-1)} \lt(|r_x-1|^2+\lt|\frac{r}{x}-1\rt|^2\rt) dx ds \\
\le & \int_0^t(1+s)^{- \frac{\ga-\ta}{\ga-\ta-\oa(\ga-1)}} \int x^2 \bar{\rho}^{\ta} \lt(|r_x-1|^2+\lt|\frac{r}{x}-1\rt|^2\rt) dx ds \\
&+ \int_0^t \int x^2 \bar{\rho}^{\ga} \lt(|r_x-1|^2+\lt|\frac{r}{x}-1\rt|^2\rt) dx ds.
\end{split}
\end{equation}
Let $\chi$ be a smooth cutoff function satisfying
\begin{align}\chi=0 \ \ {\rm on} \  \ [0,1-2\da], \ \  \chi=1 \ \ {\rm on} \ \  [1-\da,1], \ \ {\rm and} \ \  0\le \chi\le 1 \ \ {\rm on} \ \ [0,1] \label{cutoff}
\end{align}
for a small constant  $\da\in (0,1/4]$ to be determined later.
It follows from the Cauchy-Schwarz inequality  that
 \begin{equation}\label{newversion1-hh1}\begin{split}
&  \int \chi  \bar{\rho}^\ta  \lt[x\lt(|r_x-1|+\lt|\frac{r}{x}-1\rt|\rt)|v|+\lt|r-x\rt||x v_x|\rt]dx \\
 \le & \da^{1/4} (1+t) \int \chi \bar{\rho}^{\ta-(\ga-1)} v^2 dx \\
&+C \da^{-1/4} (1+t)^{-1}\int \chi x^2 \bar{\rho}^{\ta+\ga-1} \lt(|r_x-1|^2+\lt|\frac{r}{x}-1\rt|^2\rt) dx \\
&+  \da^{1/4} (1+t) \int \chi \bar{\rho}^{\ta}x^2 v_x^2 dx + C \da^{-1/4} (1+t)^{-1} \int \chi \bar{\rho}^{\ta} \lt|r-x\rt|^2 dx.
\end{split}
\end{equation}
In view of the Hardy inequality \ef{hd2} and \ef{phy}, one gets that
\begin{align*}
  \int \chi \bar{\rho}^{\ta-(\ga-1)} v^2 dx \le   \int_{1/2}^1  \bar{\rho}^{\ta-(\ga-1)} v^2 dx \le
C   \int_{1/2}^1 \bar{\rho}^{\ta+ (\ga-1)} \lt( v^2 + v_x^2\rt) dx
\end{align*}
and
\begin{align*}
\int \chi \bar{\rho}^{\ta} \lt|r-x\rt|^2 dx \le &
\int_{1-2\da}^1   \bar{\rho}^{\ta-(\ga-1)} \bar{\rho}^{ \ga-1 } \lt|r-x\rt|^2 dx \le C \da  \int_{1/2}^1   \bar{\rho}^{\ta-(\ga-1)}  \lt|r-x\rt|^2 dx   \\
\le& C \da  \int_{1/2}^1   \bar{\rho}^{\ta+(\ga-1)} \lt(|r-x|^2+ (r_x-1)^2\rt) dx,
\end{align*}
which, together with \ef{newversion1-hh1}, imply that
 \begin{equation}\label{newversion1}\begin{split}
&  \int \chi  \bar{\rho}^\ta  \lt[x\lt(|r_x-1|+\lt|\frac{r}{x}-1\rt|\rt)|v|+\lt|r-x\rt||x v_x|\rt]dx \\
 \le &   C \da^{1/4} (1+t) \int\bar{\rho}^{\ta}  \lt( v^2 + x^2v_x^2\rt) dx\\
 &+C \da^{1/4} (1+t)^{-1}\int  x^2 \bar{\rho}^{\ta+\frac{1}{2}(\ga-1)} \lt(|r_x-1|^2+\lt|\frac{r}{x}-1\rt|^2\rt) dx.
\end{split}
\end{equation}
Integrating as using \ef{hen} with $\oa=1/2$, give
 \begin{equation}\label{key}\begin{split}
 &\int_0^t\int \chi  \bar{\rho}^\ta  \lt[x\lt(|r_x-1|+\lt|\frac{r}{x}-1\rt|\rt)|v|+\lt|r-x\rt||x v_x|\rt]dxds \\
\le  &  C\da^{1/4} \int_0^t (1+s) \int\bar{\rho}^{\ta} \lt( v^2 +  x^2 v_x^2\rt) dxds \\
 &+C \da^{1/4}   \int_0^t(1+s)^{- \frac{\ga-\ta}{\ga-\ta- (\ga-1)/2}} \int x^2 \bar{\rho}^{\ta} \lt(|r_x-1|^2+\lt|\frac{r}{x}-1\rt|^2\rt) dxds \\
&+ C \da^{1/4} \int_0^t \int x^2 \bar{\rho}^{\ga} \lt(|r_x-1|^2+\lt|\frac{r}{x}-1\rt|^2\rt) dx ds.
\end{split}
\end{equation}
Using the Cauchy-Schwarz inequality and \ef{phy} again, one can obtain
\begin{equation}\label{th4}\begin{split}
 &\int_0^t \int (1-\chi) \bar{\rho}^\ta  \lt[x\lt(|r_x-1|+\lt|\frac{r}{x}-1\rt|\rt)|v|+\lt|r-x\rt||x v_x|\rt]dxds\\
  \le & C  \delta^{\frac{\ta-\ga}{\ga-1}-\frac{1}{4}}\int_0^t \int  \bar\rho^\ta  \left(x^2 v_{x}^2 + v^2  \right)  dxds\\
  & +   \delta^{\frac{\ga-\ta}{\ga-1}+\frac{1}{4}}\int_0^t \int_0^{1-\da}  \bar{\rho}^\ta x^2 \left[\left(\frac{r}{x}-1\right)^2  + \left(r_x-1\right)^2 \right] dx  ds \\
   \le & C  \delta^{\frac{\ta-\ga}{\ga-1}-\frac{1}{4}}\int_0^t \int  \bar\rho^\ta  \left(x^2 v_{x}^2 + v^2  \right)  dxds\\
  & +   C \da^{\frac{1}{4}}  \int_0^t \int_0^{1-\da}  x^2\bar{\rho}^\ga  \left[\left(\frac{r}{x}-1\right)^2  + \left(r_x-1\right)^2 \right] dx  ds\\
\le  & C  \delta^{\frac{\ta-\ga}{\ga-1}-\frac{1}{4}}\mathfrak{D}(0)+C \da^\frac{1}{4}  \int_0^t \int   x^2\bar{\rho}^\ga  \left[\left(\frac{r}{x}-1\right)^2  + \left(r_x-1\right)^2 \right] dx  ds,\end{split}
\end{equation}
where \ef{lem1est} has been used in the last inequality. So, the proof of \ef{newth3} is completed by choosing $\delta$ suitably small, and combining \ef{th3}, \ef{key}
and \ef{th4} together.

\hfill$\Box$


The following Lemma gives the decay estimates of weighted norms for the higher time-derivatives (than those in Lemma \ref{lem2}).

\begin{lem}\label{lem4} Let $\ta\in (0,1]$. Suppose that \ef{rx} and \ef{vx} hold. Then   for  $0\le t\le T$,
\begin{equation}\label{lem4est}\begin{split}
 & (1+t)  \int \lt\{ x^2 \bar{\rho} v_t^2 + \bar{\rho}^\ga\left(v^2  + x^2 v_x^2 \right) \rt\}(x,t)dx
    + \int_0^t (1+s) \int  \bar\rho^\ta  \left(x^2 v_{sx}^2 + v_s^2  \right)dx ds     \\
   & \le   C  \mathfrak{D}(0)+ C\int \lt\{ \bar{\rho} x^2v_t^2 + \bar{\rho}^\ga\left(v^2  + x^2 v_x^2 \right) \rt\}(x,0)dx.
\end{split}
\end{equation}
Moreover, for $ t\in[0,T]$, it holds that
\begin{equation}\label{lowest}\begin{split}
\mathfrak{D}_1 (t) & + \int_0^t \int x^2\bar{\rho}^\ga\left[\left(\frac{r}{x}-1\right)^2  + \left(r_x-1\right)^2 \right] dxds \\
& +   \int_0^t (1+s) \int  \bar\rho^\ta  \left(x^2 v_{x}^2 + v^2 + x^2 v_{sx}^2 + v_s^2\right)dx ds \le C \mathfrak{D}_1(0),
\end{split}
\end{equation}
where
\begin{equation*}\label{}\begin{split}
  \mathfrak{D}_1(t)   =  &  \int x^2 \bar\rho^{\ta}\lt[ \lt(r_x-1\rt)^2 +\lt(\frac{r}{x}-1\rt)^2 \rt] (x,t)dx   + (1+t)  \int x^2  \lt\{ \bar{\rho} (v^2 + v_t^2) \rt.\\
& \lt. + \bar{\rho}^\ga\left[\left(\frac{r}{x}-1\right)^2  +  \left(r_x-1\right)^2 + \lt(\frac{v}{x}\rt)^2  +   v_x^2 \right] \rt\}(x,t)dx
   .
\end{split}
\end{equation*}
\end{lem}

\noindent{\em Proof}. It yields from $\int  \lt(r^2\eqref{nsp}\rt)_t v_t dx$ that
\begin{equation}\label{vt2}\begin{split}
&  \frac{d}{dt}\int  \Phi  (x,t)dx +  \int\lt\{\nu\lt[ \lt(\frac{x^2}{r^2}\frac{\bar\rho}{ r_x}\rt)^\theta \frac{(r^2 v)_x}{r^2 r_x}\rt]_{t}(r^2 v_t)_x -4\nu_1\lt(\frac{x^2}{r^2}\frac{\bar\rho}{ r_x}\rt)^\theta\lt[(rv)_tv_t\rt]_x\rt\}dx  \\
 = &-\int \lt\{  2\nu \lt(\frac{x^2}{r^2}\frac{\bar\rho}{ r_x}\rt)^\theta\frac{(r^2v)_x}{r^2 r_x}  -4\nu_1  \lt[\lt(\frac{x^2}{r^2}\frac{\bar\rho}{ r_x}\rt)^\theta\rt]_{t}\rt\}\lt(rvv_t\rt)_x dx\\
&+  \int\lt(\frac{x^2}{r^2}\frac{\bar\rho}{ r_x}\rt)^\gamma\lt[(2\gamma-1) v_x v^2 + 2(\gamma-1) v^2 v_x + \frac \gamma 2  \lt( \frac{r^2}{r_x}\rt)_t v_x^2 \rt]dx\\
&+  \int\lt[\lt(\frac{x^2}{r^2}\frac{\bar\rho}{ r_x}\rt)^\gamma \rt]_t \lt[(2\gamma-1)  r_x  v^2 + 2(\gamma-1)  r v v_x + \frac \gamma 2   \frac{r^2}{r_x}  v_x^2 \rt]dx\\
&-  \int \bar\rho^\gamma\lt[ \lt(4\frac{x^3}{r^3}-3\frac{x^4}{r^4}r_x\rt)_t v^2 + 2\lt(\frac{x^4}{r^3}\rt)_t v v_x   \rt] dx =: {\rm RHS},
\end{split}\end{equation}
where
\begin{equation*}\label{}\begin{split}
 \Phi(x,t)&=\frac 12 x^2 \bar\rho v_t^2+\lt(\frac{x^2}{r^2} \frac{\bar\rho}{ r_x}\rt)^\gamma\lt[ (2\gamma -1) r_x v^2 + 2(\gamma-1)r vv_x  + \frac \gamma 2 \frac{r^2}{r_x} v_x^2 \rt]\\
 &-\bar\rho^\gamma \lt[\lt(4\frac{x^3}{r^3}-3\frac{x^4r_x}{r^4}\rt) v^2  +2\frac{x^4}{r^3} v v_x \rt].
\end{split}
\end{equation*}
Similar to \eqref{etalower} and \eqref{etaup}, one has
\begin{equation}\label{Philower}
{\Phi}(x,t) \ge \frac{1}{2} x^2 \bar{\rho} v_t^2 + \frac{3\ga-4}{4} x^2\bar{\rho}^\ga\left[2\left(\frac{v}{x} \right)^2  + v_x^2 \right],
\end{equation}
\begin{equation}\label{Phiup}
{\Phi}(x,t) \le \frac{1}{2} x^2 \bar{\rho} v_t^2 + C(\ga)x^2\bar{\rho}^\ga \left[\left(\frac{v}{x} \right)^2  + \left(v_x \right)^2 \right].
\end{equation}
Using the a priori assumptions \ef{rx} and \ef{vx}, and Cauchy-Schwarz's
inequality, one can obtain  that
\begin{align*}
& \int\lt\{\nu\lt[ \lt(\frac{x^2}{r^2}\frac{\bar\rho}{ r_x}\rt)^\theta \frac{(r^2 v)_x}{r^2 r_x}\rt]_{t}(r^2 v_t)_x -4\nu_1\lt(\frac{x^2}{r^2}\frac{\bar\rho}{ r_x}\rt)^\theta\lt[(rv)_tv_t\rt]_x\rt\}dx\notag \notag\\
\ge & 2\sa \int\bar\rho^{\theta}\lt(x^2v_{xt}^2  + 2 v_t^2\rt)dx-C \int\bar\rho^{\theta}\lt[ \ea_0\lt(x^2v_{xt}^2  + 2 v_t^2\rt) + \sa^{-1}   \lt(x^2v_{x}^2  +v^2\rt)\rt]dx  \end{align*}
and
\begin{align*}
{\rm RHS} \le \frac{\sigma}{2}\int\bar\rho^{\theta}\lt(x^2v_{xt}^2  +v_t^2\rt)dx+C \sa^{-1} \int \bar\rho^{\theta}\lt(x^2v_{x}^2  +v^2\rt)dx , \end{align*}
where $\sigma=\min\{ {2\nu_1}/{3}, \nu_2\}$.
Here one has used the fact that $0<\theta\le 1 < \gamma.$ Therefore, it follows from \ef{vt2}  that
\begin{equation}\label{vt2'}\begin{split}
  \frac{d}{dt}\int  \Phi &(x,t)dx +\int\bar\rho^{\theta}\lt(x^2v_{xt}^2  +v_t^2\rt) dx
\le   C \int \bar\rho^{\theta}\lt(x^2v_{x}^2  +v^2\rt)  dx, \end{split}
\end{equation}
provided that $\ea_0$ is small.
This, together with \ef{lem2est'}, implies
\be\label{7261}\int  \Phi (x,t)dx +\int_0^t \int\bar\rho^{\theta}\lt(x^2v_{xt}^2  +v_t^2\rt)(x, s)dxds\le  \int  \Phi  (x,0)dx  +C\mathfrak{D}(0), \ee
which further implies \ef{lem4est}, by using \ef{vt2'}, \ef{lem2est'}, \ef{Philower}, \ef{Phiup}  and the fact  $0<\theta\le 1<\ga$. Finally, \ef{lowest} is a  consequence of \ef{lem2est'} and \ef{lem4est}.

\hfill$\Box$

In the following lemma, we show the improved regularity near the vacuum boundary and decay.

\begin{lem}\label{lem3.4}    Suppose that \ef{rx} and \ef{vx} hold. Let $\ta\in (0,1]$, and
 $\alpha$ and $\beta$ be given  respectively in \ef{alpha} and \ef{beta}.  If $\iota\in (0, (\ga-1)/4]$,  then   for  $0\le t\le T$,
\begin{equation}\label{estlem51}\begin{split}
& \mathfrak{D}_2(t)+  \int_0^t \int x^2\bar{\rho}^{\ga-\aa} \left[\left(\frac{r}{x}-1\right)^2  + \left(r_x-1\right)^2 \right] dxds\\
&+ \int_0^t (1+s)^{ \beta-1} \int x^2\bar{\rho}^{\ga} \left[\left(\frac{r}{x}-1\right)^2  + \left(r_x-1\right)^2 \right] dxds\\
&+\int_0^t (1+s)^{ \beta} \int \bar{\rho}^\ta   \left(x^2 v_{x}^2 + v^2  \right)(x,s) dx ds \le C_\iota  \lt( \mathfrak{D}_1(0)+ \left\|  r_{0x}-1  \right\|_{L^\iy}^2 \rt),
\end{split}
\end{equation}
where
\begin{equation*}\label{}\begin{split}
  \mathfrak{D}_2(t)   =       &\int  x^2  \bar\rho^{\theta-\aa} \left[\left(\frac{r}{x}-1\right)^2  + \left(r_x-1\right)^2 \right] (x,t) dx\\
&+ (1+t)^{ \beta-1}  \int  x^2 \bar\rho^{\theta}\left[\left(\frac{r}{x}-1\right)^2  + \left(r_x-1\right)^2 \right] (x,t) dx \\
&+(1+t)^{ \beta}\int x^2\lt\{ \bar{\rho} v^2 + \bar{\rho}^\ga\left[\left(\frac{r}{x}-1\right)^2  + \left(r_x-1\right)^2 \right] \rt\}(x,t)dx
   .
\end{split}
\end{equation*}
Here $C_\iota $ is a constant continuously depending on $\iota$, but not on
$t$.

\end{lem}

\noindent{\em Proof}. The proof consists of  three steps. Recall that $\alpha$ and $\beta$ are given respectively in \ef{alpha} and \ef{beta}, that is,
$$\alpha=\min\{\ga-1+\ta, \ 2(\ga-1)\}-\iota   \ \ {\rm and} \ \    \beta=1+ ({\aa -\iota}) /({\ga-\ta}).$$

{\em Step 1}. In this step, we prove that for any $\oa>0$ and $\kappa>1$,
\be\label{htstep1}\begin{split}
   &\int x^2\bar{\rho}^{\theta-\alpha} \left[\left(\frac{r}{x}-1\right)^2  + \left(r_x-1\right)^2 \right] (x,t)  dx \\
    & +  \int_0^t \int x^2\bar{\rho}^{\ga-\alpha} \left[\left(\frac{r}{x}-1\right)^2  + \left(r_x-1\right)^2 \right] dxds\\
   \le &   C \lt( \mathfrak{D}_1(0)+ \left\|  r_{0x}-1  \right\|_{L^\iy}^2 \rt)+ C \oa \int_0^t (1+s)^\kappa  \int \bar\rho^{\theta} (v^2 + x^2 v_x^2) dx  ds \\
 & +  C\oa^{-1} \int_0^t (1+s)^{-\kappa} \int x^2 \bar\rho^{\ta-\alpha} \left[\left(\frac{r}{x}-1\right)^2  + \left(r_x-1\right)^2 \right]dx  ds.
\end{split}\ee
It should be noted that $\kappa>1$, which ensures the last line in \ef{htstep1} behaves well according to the Gronwall inequality.

Multiplying \eqref{nsp} by $\int_0^x \bar\rho^{- \alpha}(y)\lt(r^3-y^3\rt)_y dy$ and integrating the resulting equation with respect to the spatial variable, one has, with the help of the integration by parts and  boundary condition \eqref{lex1}, that
 \begin{equation*}\label{}\begin{split}
 & \int     \bar{\rho}^\ga \left\{       \left[\frac{x^4}{r^4} \int_0^x \bar\rho^{- \alpha}(y)\lt(r^3-y^3\rt)_y dy\right]_x  - \left(\frac{x^2}{r^2 r_x} \right)^\ga  \bar\rho^{\ga- \alpha}  \left(r^3-x^3\right)_x   \right\} dx\\
   = &   -\int \bar\rho\left( \frac{x}{r}\right)^2  v_t \int_0^x \bar\rho^{- \alpha}(y)\lt(r^3-y^3\rt)_y dy dx
- \nu  \int  \bar{\rho}^{\ta - \alpha}  \left(\frac{x^2}{r^2 r_x}\right)^\ta \frac{ (r^2v)_x}{r^2  r_x}\left(r^3-x^3\right)_x dx\\
& + 4\nu_1 \int \bar{\rho}^\ta  \left(\frac{x^2}{r^2 r_x}\right)^\ta \lt[\frac{v}{ r}\int_0^x \bar\rho^{- \alpha}(y)\lt(r^3-y^3\rt)_y dy\rt]_x   dx   .
\end{split}
\end{equation*}
As the derivation of \ef{th2}, one can get
\be\label{beauty1}\begin{split}
   &\int x^2\bar{\rho}^{\theta-\alpha} \left[\left(\frac{r}{x}-1\right)^2  + \left(r_x-1\right)^2 \right]  dx  +  \int_0^t \int x^2\bar{\rho}^{\ga-\alpha} \left[\left(\frac{r}{x}-1\right)^2  + \left(r_x-1\right)^2 \right] dxds\\
  & \le C \int x^2\bar{\rho}^{\theta-\alpha} \left[\left(\frac{r}{x}-1\right)^2  + \left(r_x-1\right)^2 \right](x,0) dx   + C \sum_{i=1}^3 \int_0^t |L_i| ds,
\end{split}\ee
where
\begin{equation*}\label{Li}\begin{split}
L_1= &-\int \bar\rho \frac{x^2}{r^2}  v_t   \left( \int_0^x \bar\rho^{-\alpha}(r^3-y^3)_ydy \right) dx , \\
L_2=& \int     \bar{\rho}^\ga \lt(\frac{x^4}{r^4}\rt)_x \left[\bar\rho^{-\alpha} \left(r^3-x^3\right) -  \int_0^x \bar\rho^{-\alpha}(r^3-y^3)_ydy  \right] dx, \\
L_3=& 4\nu_1 \int \bar{\rho}^\ta  \left(\frac{x^2}{r^2 r_x}\right)^\ta \lt[\frac{v}{ r}\int_0^x \bar\rho^{- \alpha}(y)\lt(r^3-y^3\rt)_y dy\rt]_x   dx.
\end{split}
\end{equation*}
In the above integration by parts, we should notice that $0<\alpha<\ta+ (\ga-1)$ which ensures that $\bar\rho^{\ta-\alpha}$ is integrable on $[0, 1]$ due
to \ef{phy}, so that the above integrations by parts are legitimate.

Next, we estimate $L_1$, $L_2$ and $L_3$. For $L_1$, it follows from the Cauchy inequality that for any $\oa>0$,
\begin{align}
&|L_1|\le   C\int \bar\rho  v_t   \left| \int_0^x \bar\rho^{-\alpha}y^2(|r/y-1| +   |r_y -1|  )dy \right| dx \notag\\
\le & C \oa^{-1} \int \bar\rho^{\ta-(\ga-1)}  v_t^2 dx + \oa \int {\bar\rho}^{1-\ta+\ga}   \left| \int_0^x \bar\rho^{-\alpha}y^2 (|r/y-1| +   |r_y -1|  )dy \right|^2 dx. \label{L1-1}
\end{align}
Due to \ef{phy} and \ef{hd2}, one has
\begin{align}
&\int \bar\rho^{\ta-(\ga-1)}  v_t^2 dx \le C \int_0^{1/2} \bar\rho^{\ta}  v_t^2 dx + C \int_{1/2}^1 \bar\rho^{\ta+(\ga-1)} ( v_t^2 + v_{tx}^2 ) dx \notag \\
\le &  C \int_0^{1/2} \bar\rho^{\ta}  v_t^2 dx + C \int_{1/2}^1  \bar\rho^{\ta} ( v_t^2 + x^2 v_{tx}^2 ) dx \le C \int  \bar\rho^{\ta} ( v_t^2 + x^2 v_{tx}^2 ) dx. \label{L1-2}
\end{align}
Due to H\"{o}lder's inequality, $\ta\le  1$ and $\alpha<2(\ga-1)$ (which implies $\alpha+\ta < 2\ga -1$), and \ef{phy}, one obtains
\begin{align}
& \int {\bar\rho}^{1-\ta+\ga}   \left| \int_0^x \bar\rho^{-\alpha}y^2 (|r/y-1| +   |r_y -1|  )dy \right|^2 dx \notag\\
 \le &
 \int {\bar\rho}^{1-\ta+\ga}    \lt( \int_0^x \bar\rho^{-\ga -\alpha} y^2  dy \rt) \lt( \int_0^1 \bar\rho^{\ga-\alpha} y^2 (|r/y-1| +   |r_y -1|  )dy \rt) dx \notag\\
 \le  & C   \int_0^1 \bar\rho^{\ga-\alpha} y^2 (|r/y-1| +   |r_y -1|  )dy. \label{L1-3}
\end{align}
It yields from \ef{L1-1}-\ef{L1-3} that for any $\oa>0$,
\begin{align}
\int_0^t |L_1|  ds \le & C \oa^{-1}\int_0^t  \int  \bar\rho^{\ta} ( v_s^2 + x^2 v_{sx}^2 ) dx ds \notag\\
& + C\oa \int_0^t \int x^2\bar{\rho}^{\ga-\alpha} \left[\left(\frac{r}{x}-1\right)^2  + \left(r_x-1\right)^2 \right] dxds . \label{newhL1}
\end{align}
Rewrite $L_2$ as
\begin{equation*}\label{}\begin{split}
L_2 =& \int_0^{1/2} \bar{\rho}^{\ga}\lt(\frac{x^4}{r^4}\rt)_x \int_0^x \lt(\bar\rho^{-\alpha}\rt)_y(r^3-y^3)dydx \\
 &+ \int_{1/2}^1 \bar{\rho}^{\ga}\lt(\frac{x^4}{r^4}\rt)_x\lt[\bar{\rho}^{-\alpha}(r^3-x^3)
 - \int_0^x \bar\rho^{-\alpha}(r^3-y^3)_ydy\rt]dx=:L_{21}+L_{22}.
\end{split}
\end{equation*}
It follows from  that \ef{rhox}, and the Cauchy and H\"{o}lder inequalities that
\begin{equation}\label{l21}\begin{split}
& |L_{21}| \le  C \int_0^{1/2} x^{-1} \lt(|r_x-1|+\lt|\frac{r}{x}-1\rt|\rt) \lt|\int_0^x y^3 |r-y| dy\rt| dx\\
\le &  C \int_0^{1/2} x^2 \lt(|r_x-1|^2+\lt|\frac{r}{x}-1\rt|^2\rt) dx + C \int_0^{1/2} x^{-4} \lt(\int_0^x y^6 dy \rt) \lt( \int_0^x |r-y|^2 dy \rt)  dx \\
\le & C \int_0^{1/2} x^2 \lt(|r_x-1|^2+\lt|\frac{r}{x}-1\rt|^2\rt) dx \le C \int_0^{1/2} x^2 \bar\rho^\ga \lt(|r_x-1|^2+\lt|\frac{r}{x}-1\rt|^2\rt) dx,
\end{split}
\end{equation}
and for any $\oa>0$,
\begin{equation}\label{l22}\begin{split}
&|L_{22}| \le   C\int_{1/2}^1 \bar{\rho}^{\ga-\alpha}\lt(|r_x-1|+\lt|\frac{r}{x}-1\rt|\rt)|r-x|dx \\
 & \quad   +C \int_{1/2}^1  \bar{\rho}^{\ga}\lt(|r_x-1|+\lt|\frac{r}{x}-1\rt|\rt)\lt(\int_0^x \bar\rho^{-\alpha}y^2\lt(|r_y-1|+\lt|\frac{r}{y}-1\rt|\rt)dy \rt) dx\\
 \le &   \oa \int _{1/2}^1 \bar{\rho}^{\ga-\alpha} \left[\left(\frac{r}{x}-1\right)^2  + \left(r_x-1\right)^2 \right] dx
 + C \oa^{-1} \int _{1/2}^1 \bar{\rho}^{\ga-\alpha} (r-x)^2 dx \\
&  + C \oa^{-1} \int_{1/2}^1  \bar{\rho}^{\ga +\alpha } \lt(\int_0^x \bar\rho^{-\ga-2\alpha}y^2 dy\rt) \lt( \int_0^1 \bar\rho^{\ga }y^2\lt(|r_y-1|^2+\lt|\frac{r}{y}-1\rt|^2\rt)dy \rt) dx.
\end{split}
\end{equation}
By virtue of  \ef{phy},  \ef{hd2}  and   $\alpha < 2(\ga-1)$,  one has
\begin{equation*}\label{subcase1}\begin{split}
\int_{1/2}^1 \bar\rho^{\ga-\alpha} (r-x)^2 dx \le & C \int_{1/2}^1 \bar\rho^{\ga-\alpha + 2(\ga-1)} \lt[(r-x)^2 +(r_x-1)^2\rt] \\
\le & C \int_{1/2}^1 \bar\rho^{\ga} \lt[(r-x)^2 +x^2 (r_x-1)^2\rt]
\end{split}
\end{equation*}
and
$$\int_{1/2}^1  \bar{\rho}^{\ga +\alpha } \lt(\int_0^x \bar\rho^{-\ga-2\alpha}y^2 dy\rt)  dx \le C.$$
It then yields from \ef{l21} and \ef{l22}   that for any $\oa>0$,
\begin{equation}\label{14l2}\begin{split}
\int_0^t  |L_2|ds
\le & C\oa \int_0^t \int x^2\bar{\rho}^{\ga-\alpha} \left[\left(\frac{r}{x}-1\right)^2  + \left(r_x-1\right)^2 \right] dxds \\
& + C\oa^{-1} \int_0^t \int x^2\bar{\rho}^{\ga} \left[\left(\frac{r}{x}-1\right)^2  + \left(r_x-1\right)^2 \right] dxds.
\end{split}
\end{equation}
For $L_3$, it holds that
\begin{equation}\label{l3}\begin{split}
|L_3|\le  &C \int  \bar\rho^{\theta-\alpha} \lt| \frac{v}{r} \rt| \lt| (r^3-x^3)_x  \rt|dx +  C \int \bar\rho^{\theta}  \left |\lt(\frac{v}{r}\rt)_x\rt |  \lt| \int_0^x \bar\rho^{-\alpha} (r^3-y^3)_y  dy \rt| dx   \\
\le & C \int  \bar\rho^{\theta-\alpha} |v|   \lt( x|r_x-1| + |r-x|  \rt)dx    \\
 & + C \int \bar\rho^{\theta} x^{-2} (x |v_x| + |v|)   \lt| \int_0^x \bar\rho^{-\alpha} (r^3-y^3)_y  dy \rt| dx
 =: L_{31}+L_{32}.
\end{split}
\end{equation}
It follows from \ef{hd1}, \ef{hd2}, \ef{phy},  $\ta-\alpha>1-\ga$ and   $\alpha < 2(\ga-1)$ that
\begin{equation}\label{5/10-1}\begin{split}
 \int \bar\rho^{\theta-\alpha} v^2 dx  \le C \int_0^{1/2}  v^2 dx + \int_{1/2}^1 \bar\rho^{\theta-\alpha} v^2 dx
 \le
 C\int_0^{1/2} x^2 ( v^2 + v_x^2 ) dx  \\
 + C\int_{1/2}^1 \bar\rho^{\theta-\alpha+ 2(\ga-1)} (v^2 + v_x^2) dx \le C\int \bar\rho^{\theta} (v^2 + x^2 v_x^2) dx ,
\end{split}
\end{equation}
which, together with the Cauchy inequality, gives that for any $\oa>0$,
\begin{equation}\label{l31}\begin{split}
L_{31}\le &   C \oa (1+t)^\kappa  \int \bar\rho^{\theta} (v^2 + x^2 v_x^2) dx \\
& +  C\oa^{-1} (1+t)^{-\kappa} \int x^2 \bar\rho^{\ta-\alpha} \left[\left(\frac{r}{x}-1\right)^2  + \left(r_x-1\right)^2 \right]dx  .
\end{split}
\end{equation}
Clearly, $L_{32}$ can be bounded by
\begin{equation*}\label{}\begin{split}
& L_{32} \le   C \int_0^{1/2}  x^{-2} (x |v_x| + |v|)   \lt|\bar\rho^{-\alpha} (r^3 -x^3) -  \int_0^x (\bar\rho^{-\alpha})_y (r^3-y^3)   dy \rt| dx \\
& + C\int_{1/2}^1 \bar\rho^{\theta}  (x |v_x| + |v|)   \lt| \int_0^x \bar\rho^{-\alpha} y^2 \lt(|r_y-1|+\lt|\frac{r}{y}-1\rt|\rt)dy   \rt| dx
=:   L_{321}+L_{322}.
\end{split}
\end{equation*}
For $L_{321}$, it follows from \ef{rhox} and the Cauchy and H\"{o}lder inequalities that
 \begin{equation*}\label{l321}\begin{split}
& L_{321}\le     C \int_0^{1/2}   (x |v_x| + |v|) |r-x| dx +  C \int_0^{1/2}   x^{-2} (x |v_x| + |v|)
   \lt| \int_0^x  y^3 |r-y|  dy \rt| dx  \\
   \le &  C \int_0^{1/2} \lt[\bar\rho^\ta (x^2 v_x^2 + v^2) + \bar\rho^\ga (r-x)^2 \rt]dx
   + C \int_0^{1/2}   x^{-4} \lt(\int_0^x  y^6 dy \rt) \lt( \int_0^x (r-y)^2  dy\rt)  dx \\
   \le &  C \int_0^{1/2} \lt[\bar\rho^\ta (x^2 v_x^2 + v^2) + \bar\rho^\ga (r-x)^2 \rt]dx.
\end{split}
\end{equation*}
For $L_{322}$, it follows from  the Cauchy and H\"{o}lder inequalities that for any $\oa>0$,
\begin{equation*}\label{}\begin{split}
L_{322} \le & \oa  (1+t)^{\kappa} \int_{1/2}^1 \bar\rho^{\theta}  (x^2 v_x^2 + v^2) dx \\
&+ C \oa^{-1} (1+t)^{-\kappa}  \int_0^1 \bar\rho^{\ta-\alpha} y^2 \lt(|r_y-1|^2+\lt| {r}/{y}-1\rt|^2\rt)dy
.
\end{split}
\end{equation*}
Here one has used the fact that $\alpha < 2(\ga-1)$, which implies
$$ \int_{1/2}^1 \bar\rho^{\theta}   \lt( \int_0^x \bar\rho^{-\ta-\alpha} y^2 dy \rt) dx \le C.$$
So, it holds that for any $\oa>0$,
\begin{equation*}\label{}\begin{split}
L_{32}\le &    \oa (1+t)^\kappa  \int \bar\rho^{\theta} (v^2 + x^2 v_x^2) dx + C \int  \lt[\bar\rho^\ta (x^2 v_x^2 + v^2) + \bar\rho^\ga (r-x)^2 \rt]dx \\
& +  C\oa^{-1} (1+t)^{-\kappa} \int x^2 \bar\rho^{\ta-\alpha} \left[\left(\frac{r}{x}-1\right)^2  + \left(r_x-1\right)^2 \right]dx  .
\end{split}
\end{equation*}
This, together with \ef{l3} and \ef{l31}, implies
\begin{equation}\label{newl3}\begin{split}
\int_0^t |L_{3}| ds \le &   C \oa \int_0^t (1+s)^\kappa  \int \bar\rho^{\theta} (v^2 + x^2 v_x^2) dx  ds \\
&+ C \int_0^t \int  \lt[\bar\rho^\ta (x^2 v_x^2 + v^2) + \bar\rho^\ga (r-x)^2 \rt]dx ds \\
& +  C\oa^{-1} \int_0^t (1+s)^{-\kappa} \int x^2 \bar\rho^{\ta-\alpha} \left[\left(\frac{r}{x}-1\right)^2  + \left(r_x-1\right)^2 \right]dx  ds .
\end{split}
\end{equation}

The estimate \ef{htstep1} follows from \ef{newhL1}, \ef{14l2}, \ef{newl3}, \ef{lowest}, and the fact that $r(0,t)=0$ (which implies $\|r-x\|_{L^\iy} \le \|r_x -1\|_{L^\iy} $).

{\em Step 2}.  In this step, we prove that for any $\kappa>1$,
\begin{align}
\mathscr{F}_1^\kappa(t) & \le C\mathfrak{D}_1(0)+
C \int_0^t (1+s)^{-1+\lt(\kappa-1-\frac{\alpha}{\ga-\ta}\rt)}
\int x^2 \bar\rho^{\theta-\aa} \left[\left(\frac{r}{x}-1\right)^2  + \left(r_x-1\right)^2 \right]  dxds \notag \\
& +C  \int_0^t(1+s)^{ - \frac{\ga-\ta}{\ga-\ta-(\ga-1)/2}} (1+s)^{ \kappa-1} \int x^2 \bar{\rho}^{\ta} \lt(|r_x-1|^2+\lt|\frac{r}{x}-1\rt|^2\rt) dx ds \notag\\
& + C (1+t)^{\kappa-1} \int x^2 \bar\rho v^2(x,t) dx + C \int_0^t (1+s)^{\kappa-1} \int \bar\rho^\ta (x^2 v_x^2+ v^2 )dxds, \label{keyconclusion}
\end{align}
where
\begin{align*}
\mathscr{F}_1^\kappa(t)= &(1+t)^{ \kappa}\int x^2\lt\{ \bar{\rho} v^2 + \bar{\rho}^\ga\left[\left(\frac{r}{x}-1\right)^2  + \left(r_x-1\right)^2 \right] \rt\}(x,t)dx\notag\\
&+ (1+t)^{ \kappa-1}  \int  x^2 \bar\rho^{\theta}\left[\left(\frac{r}{x}-1\right)^2  + \left(r_x-1\right)^2 \right] (x,t) dx\notag\\
&+\int_0^t (1+s)^{ \kappa} \int \bar{\rho}^\ta   \left(x^2 v_{x}^2 + v^2  \right)  dx ds
\notag\\
&+  \int_0^t (1+s)^{ \kappa-1} \int x^2\bar{\rho}^{\ga} \left[\left(\frac{r}{x}-1\right)^2  + \left(r_x-1\right)^2 \right] dxds.
\end{align*}

It follows from a suitable combination of  $\int_0^t (1+s)^{ \kappa-1}  \ef{th1}ds$ with \ef{5-0-1} for $\ell= \kappa$ that
\begin{align}\label{cc2}
\mathscr{F}_1^\kappa(t) & \le C\mathfrak{D}_1(0)+
 C \int_0^t(1+s)^{ \kappa-2}\int x^2\bar\rho^{\theta} \left[\left(\frac{r}{x}-1\right)^2  + \left(r_x-1\right)^2 \right]  dxds\notag\\
&+C \int_0^t (1+s)^{ \kappa-1} \int\lt[ \bar{\rho}^\ta x\lt(|r_x-1|+\lt|\frac{r}{x}-1\rt|\rt)|v|+  \bar{\rho}^\ta \lt|r-x\rt||xv_x|\rt]dxds\notag\\
& + C (1+t)^{\kappa-1} \int x^2 \bar\rho v^2(x,t) dx + C \int_0^t (1+s)^{\kappa-1} \int \bar\rho^\ta v^2 dxds,
\end{align}
where \ef{cc1} and \ef{lowest} have been used. Each term on the right-hand side of \ef{cc2} can be estimated as follows. It follows from the Young inequality  that for any $\oa>0$,
\be\label{cc3}\begin{split}
&\int_0^t(1+s)^{ \kappa-2}\int x^2 \bar\rho^{\theta}\left[\left(\frac{r}{x}-1\right)^2  + \left(r_x-1\right)^2 \right]  dxds\\
 \le & \oa \int_0^t (1+s)^{\kappa-1} \int x^2 \bar\rho^{\ga} \left[\left(\frac{r}{x}-1\right)^2  + \left(r_x-1\right)^2 \right] dxds
\\
&+
C_\oa\int_0^t (1+s)^{-1+\lt(\kappa-1-\frac{\alpha}{\ga-\ta}\rt)}
\int x^2 \bar\rho^{\theta-\aa} \left[\left(\frac{r}{x}-1\right)^2  + \left(r_x-1\right)^2 \right]  dxds,
\end{split}\ee
due to
\begin{align*}
 (1+s)^{ \kappa-2}  \bar\rho^{\theta} = \lt((1+s)^{ \frac{\ga-\ta}{\ga-\ta+\aa} (\kappa-1)-1}  \bar\rho^{\frac{\ga-\ta}{\ga-\ta+\aa}(\theta-\aa)} \rt) \lt( (1+s)^{ \frac{\aa}{\ga-\ta+\aa} ( \kappa-1)}  \bar\rho^{\frac{\aa}{\ga-\ta+\aa}\ga}   \rt).
\end{align*}
Let $\chi$ be a smooth cutoff function satisfying \ef{cutoff}. Similar to \ef{newversion1} and \ef{hen}, one  obtains
\begin{equation}\label{newversion2}\begin{split}
&  (1+t)^{\kappa-1}\int \chi  \bar{\rho}^\ta  \lt[x\lt(|r_x-1|+\lt|\frac{r}{x}-1\rt|\rt)|v|+\lt|r-x\rt||x v_x|\rt]dx \\
 \le &  C \da^{1/4} (1+t)^{\kappa} \int\bar{\rho}^{\ta}  \lt( v^2 + x^2v_x^2\rt) dx\\
 &+C \da^{1/4} (1+t)^{\kappa-2}\int  x^2 \bar{\rho}^{\ta+ (\ga-1)/2} \lt(|r_x-1|^2+\lt|\frac{r}{x}-1\rt|^2\rt) dx
\end{split}
\end{equation}
and
\begin{equation}\label{newversion3}\begin{split}
  & (1+t)^{\kappa-2}\int  x^2 \bar{\rho}^{\ta+ (\ga-1)/2} \lt(|r_x-1|^2+\lt|\frac{r}{x}-1\rt|^2\rt) dx\\
 \le &  (1+t)^{\kappa-1- \frac{\ga-\ta}{\ga-\ta-(\ga-1)/2}}\int x^2 \bar{\rho}^{\ta} \lt(|r_x-1|^2+\lt|\frac{r}{x}-1\rt|^2\rt) dx\\
& + (1+t)^{\kappa-1}\int x^2 \bar{\rho}^{\ga} \lt(|r_x-1|^2+\lt|\frac{r}{x}-1\rt|^2\rt) dx.
\end{split}
\end{equation}
It thus yields from \ef{newversion2} and \ef{newversion3} that
\begin{equation}\label{keyx}\begin{split}
& \int_0^t(1+s)^{\kappa-1} \int \chi \bar{\rho}^\ta  \lt[x\lt(|r_x-1|+\lt|\frac{r}{x}-1\rt|\rt)|v|+\lt|r-x\rt||x v_x|\rt]dxds \\
 \le  &  C\da^{1/4} \int_0^t (1+s)^{\kappa} \int\bar{\rho}^{\ta} \lt( v^2 + x^2v_x^2\rt) dxds \\
 &+C \da^{1/4}   \int_0^t(1+s)^{ \kappa-1- \frac{\ga-\ta}{\ga-\ta-(\ga-1)/2}} \int x^2 \bar{\rho}^{\ta} \lt(|r_x-1|^2+\lt|\frac{r}{x}-1\rt|^2\rt) dx ds\\
&+ C \da^{1/4} \int_0^t (1+s)^{\kappa-1}\int x^2 \bar{\rho}^{\ga} \lt(|r_x-1|^2+\lt|\frac{r}{x}-1\rt|^2\rt) dx ds.
\end{split}
\end{equation}
Similar to \ef{th4}, one can get
\begin{equation}\label{th41}\begin{split}
 &\int_0^t (1+s)^{\kappa-1} \int (1-\chi) \bar{\rho}^\ta  \lt[x\lt(|r_x-1|+\lt|\frac{r}{x}-1\rt|\rt)|v|+\lt|r-x\rt||x v_x|\rt]dxds\\
      \le & C  \delta^{\frac{\ta-\ga}{\ga-1}-\frac{1}{4}}\int_0^t (1+s)^{\kappa-1} \int  \bar\rho^\ta  \left(x^2 v_{x}^2 + v^2  \right)  dxds\\
  & +   C \da^{\frac{1}{4}}  \int_0^t (1+s)^{\kappa-1} \int   x^2\bar{\rho}^\ga  \left[\left(\frac{r}{x}-1\right)^2  + \left(r_x-1\right)^2 \right] dx  ds.\\
  \end{split}
\end{equation}
Now, the proof of \ef{keyconclusion} is completed by choosing $\oa$ and $\da$ suitably small, and combining \ef{cc2}, \ef{cc3}, \ef{keyx} and \ef{th41}.

{\em Step 3}. In this step, we prove that
 \begin{align}
 \mathscr{F}_2^\beta(t)   \le & C \lt( \mathfrak{D}_1(0)+ \left\|  r_{0x}-1  \right\|_{L^\iy}^2 \rt) , \label{newhtstep3}
\end{align}
where $\beta$ is defined by \ef{beta}, and
\begin{align*}
\mathscr{F}_{2}^\kappa(t)= &  \mathscr{F}_1^\kappa(t) +
 \int x^2\bar{\rho}^{\theta-\alpha} \left[\left(\frac{r}{x}-1\right)^2  + \left(r_x-1\right)^2 \right] (x,t)  dx  \\
    & +  \int_0^t \int x^2\bar{\rho}^{\ga-\alpha} \left[\left(\frac{r}{x}-1\right)^2  + \left(r_x-1\right)^2 \right] dxds  \ \ {\rm for} \ \ {\rm any} \ \  \kappa>1.
\end{align*}

It follows from \ef{htstep1} and \ef{keyconclusion} that for any $\kappa>1$,
\begin{align}
 \mathscr{F}_{2}^\kappa (t)   \le & C \lt( \mathfrak{D}_1(0)+ \left\|  r_{0x}-1  \right\|_{L^\iy}^2 \rt)  + C \int_0^t \lt[ (1+s)^{-1+\lt(\kappa-1-\frac{\alpha}{\ga-\ta}\rt)} \rt.\notag\\
& \lt. +(1+s)^{-\kappa} \rt]
\int x^2 \bar\rho^{\theta-\aa} \left[\left(\frac{r}{x}-1\right)^2  + \left(r_x-1\right)^2 \right]  dxds \notag \\
& +C  \int_0^t(1+s)^{ - \frac{\ga-\ta}{\ga-\ta-(\ga-1)/2}} (1+s)^{ \kappa-1} \int x^2 \bar{\rho}^{\ta} \lt(|r_x-1|^2+\lt|\frac{r}{x}-1\rt|^2\rt) dx ds \notag\\
& + C (1+t)^{\kappa-1} \int x^2 \bar\rho v^2(x,t) dx + C \int_0^t (1+s)^{\kappa-1} \int \bar\rho^\ta (x^2 v_x^2+ v^2 )dxds  . \label{htstep3}
\end{align}
When $\beta \le 2$, \ef{newhtstep3} follows from
follows from \ef{htstep3} with $\kappa = \beta$, \ef{lowest},  $\ta\le  1$, and the Grownwall inequality. When $\beta>2$, we have
$$1-\aa/(\ga-\ta)<-\iota/(\ga-\ta).$$
Then, it follows from \ef{htstep3} with $\kappa = 2$, \ef{lowest},  $\ta\le   1$, and the Grownwall inequality that
\begin{align}
 \mathscr{F}_{2}^2 (t)   \le & C \lt( \mathfrak{D}_1(0)+ \left\|  r_{0x}-1  \right\|_{L^\iy}^2 \rt). \label{hala}
\end{align}
Due to $\ta\le 1$, one has
$$\beta=1+\frac{\alpha-\iota}{\ga-\ta} < 1 + \frac{\alpha}{\ga-1} <3 .$$
So, \ef{newhtstep3} follows from \ef{htstep3} with $\kappa = \beta$, \ef{hala},  $\ta\le  1$, and the Grownwall inequality.

\hfill$\Box$

As an immediate consequence of \ef{estlem51}, we have the following Lemma.

\begin{lem}\label{lem3.5} Suppose that \ef{rx} and \ef{vx} hold. Let $\ta\in (0,\ga/2]$, and
 $\alpha$, $\beta$  and $\varsigma$ be given  respectively in \ef{alpha}, \ef{beta} and \ef{varsa}.  If $\iota\in (0, (2\ga-2-\ta)/8 ]$, then   for  $0\le t\le T$,
 \begin{equation}\label{3.63}\begin{split}
& (1+t)^{\frac{2\ga-2+\aa-\theta}{\ga+\aa-\ta}\beta} \int (r(x, t)-x)^2dx + (1+t)^{\frac{ \aa-\theta}{\ga+\aa-\ta}\beta} \int x^2 (r_x(x, t)-1)^2dx \\
&+ (1+t)^{  \beta}\int \lt(x^2\bar\rho v_t^2+\bar\rho^\ga(x^2 v_x^2+v^2)\rt)(x, t)dx\\
&    +\int_0^t (1+s)^{ \beta}\int\bar\rho^{\theta}(x^2v_{xs}^2+v_s^2) dxds
 \le C_\iota  \lt( \mathfrak{D}_1(0)+ \left\|  r_{0x}-1  \right\|_{L^\iy}^2 \rt),
\end{split}\end{equation}
\begin{equation}\label{3.64}\begin{split}
&(1+t)^\beta \int \lt(v^2 + x^2 \bar\rho^{\theta} v_{x}^2 \rt) (x,t) dx
 + (1+t)^{\frac{\beta+\varsigma}{2}  } \int\bar\rho^{\theta- \aa/2}(x^2v_{x}^2+v^2) (x,t) dx \\
&  \le  C_\iota  \lt( \mathfrak{D}_1(0)+ \left\|  r_{0x}-1  \right\|_{L^\iy}^2 + \left\| v_x(\cdot, 0) \right\|_{L^\iy}^2 \rt).
\end{split}\end{equation}

\end{lem}
{\em Proof}. The proof consists of three steps.

{\em Step 1}. In this step, we prove that
\begin{equation}\label{estlem512}\begin{split}
(1+t)^{\frac{2\ga-2+\aa-\theta}{\ga+\aa-\ta}\beta} \int (r(x, t)-x)^2dx\le C_\iota  \lt( \mathfrak{D}_1(0)+ \left\|  r_{0x}-1  \right\|_{L^\iy}^2 \rt),
\end{split}\end{equation}
\begin{equation}\label{5/10-2}\begin{split}
(1+t)^{\frac{ \aa-\theta}{\ga+\aa-\ta}\beta} \int x^2 (r_x(x, t)-1)^2dx\le C_\iota  \lt( \mathfrak{D}_1(0)+ \left\|  r_{0x}-1  \right\|_{L^\iy}^2 \rt),
\end{split}\end{equation}
\begin{equation}\label{estlem51'}\begin{split}
&(1+t)^{  \beta}\int \lt(x^2\bar\rho v_t^2+\bar\rho^\ga(x^2 v_x^2+v^2)\rt)(x, t)dx\\
& \quad\quad  +\int_0^t (1+s)^{ \beta}\int\bar\rho^{\theta}(x^2v_{xs}^2+v_s^2) dxds  \le C_\iota  \lt( \mathfrak{D}_1(0)+ \left\|  r_{0x}-1  \right\|_{L^\iy}^2 \rt).
\end{split}\end{equation}
It follows from Hardy's inequalities \ef{hd1} and \ef{hd2}, \ef{phy}, $\ta<\alpha$ and the H\"{o}lder inequality that
\begin{equation*}\begin{split}
   \int(r-x)^2 dx \le & C \int x^2 \bar\rho^{2(\ga-1)}\left((r-x)^2+( r_x-1)^2 \right)dx \\
   \le &  C \lt(\int x^2 \bar\rho^{\theta-\aa}\left((r-x)^2+( r_x-1)^2  \right)dx\rt)^{\frac{2-\ga}{\ga+\aa-\ta}} \\
   &\times
 \lt( \int x^2 \bar\rho^{\ga}\left((r-x)^2+( r_x-1)^2 \right)dx \rt)^{\frac{2\ga-2+\aa-\ta}{\ga+\aa-\ta}}.
  \end{split}
\end{equation*}
Thus, \ef{estlem512} follows from \ef{estlem51}. Clearly, \ef{5/10-2} follows from \ef{estlem51} and
\begin{equation*}\begin{split}
   \int   x^2 ( r_x-1)^2  dx
   \le &  C \lt(\int  x^2  \bar\rho^{\theta-\aa}  ( r_x-1)^2  dx\rt)^{\frac{ \ga}{\ga+\aa-\ta}}
 \lt( \int   x^2 \bar\rho^{\ga}    ( r_x-1)^2   dx \rt)^{\frac{\aa-\ta}{\ga+\aa-\ta}}.
  \end{split}
\end{equation*}
With \ef{estlem51}, one may easily obtain \ef{estlem51'} by virtue of \ef{vt2'}, \ef{Philower}, \ef{Phiup} and \ef{lowest}.

{\em Step 2}. In this step, we prove that
\begin{equation}\label{511-1}\begin{split}
&(1+t)^{ \varsigma} \int \bar\rho^{\ga-\aa} \lt[ (x^2(r_x-1)^2+ (r-x)^2 )\rt](x, t)dx\\
& +\int_0^t (1+s)^{ \varsigma }\int\bar\rho^{\theta-\aa}(x^2v_{x }^2+v^2) dxds  \le C_\iota  \lt( \mathfrak{D}_1(0)+ \left\|  r_{0x}-1  \right\|_{L^\iy}^2 \rt).
\end{split}\end{equation}
Multiplying equation \eqref{nsp} by $\int_0^x \bar\rho^{-\aa}(y)(r^2 v)_ydy$ and integrating the product with respect to spatial variable,  one obtains, following the derivation of \ef{5-0}, that
\begin{equation}\label{beauty}\begin{split}
&(1+t)^{\varsigma}\int x^2\bar{\rho}^{\ga-\aa}\left[\left(\frac{r}{x}-1\right)^2  + \left(r_x-1\right)^2 \right] (x,t) dx
+ \int_0^t  (1+s)^{\varsigma}  \int \bar\rho^{\theta-\aa} \left(x^2 v_x^2+ v^2  \right)dxds\\
& \le C \left\|  r_{0x}-1  \right\|_{L^\iy}^2 +  C \sum_{i=1}^3 \int_0^t (1+s)^{\varsigma}| K_i |ds
 +  C\int_0^t  \int x^2\bar{\rho}^{\ga-\aa}\left[\left(\frac{r}{x}-1\right)^2  + \left(r_x-1\right)^2 \right] dxds,
\end{split}
\end{equation}
where
\begin{equation*}\label{Ki}\begin{split}
K_1  = & \int \bar\rho\frac{x^2}{r^2}v_t \lt(\int_0^x \bar\rho^{-\aa}(r^2 v)_ydy \rt)dx,\\
K_2= & \int \bar{\rho}^{\ga}\lt(\frac{x^4}{r^4}\rt)_x\lt[\bar{\rho}^{-\aa}{r^2}v
 - \int_0^x \bar\rho^{-\aa}(r^2 v)_ydy\rt]dx,\\
K_3= & 4 \nu_1 \int \lt(\frac{x^2}{r^2}\frac{\bar\rho}{ r_x}\rt)^\theta\left(\frac{v}{r}\right)_x \lt[ \lt(\int_0^x \bar\rho^{-\aa}(r^2 v)_ydy \rt)-\bar\rho^{-\aa} r^2 v \rt]dx.
\end{split}
\end{equation*}
The estimate for $K_i$ can be obtained in a similar way as the derivation of \ef{newhL1}, \ef{14l2} and \ef{newl3}. Clearly,
\begin{equation}\label{estK1}\begin{split}
\int_0^t (1+s)^{\varsigma} | K_1 |   \le & C \oa^{-1} \int_0^t (1+s)^{\varsigma} \int \bar\rho^\ta v_s^2 dxds\\& + C\oa \int_0^t (1+s)^{\varsigma} \int \bar\rho^{\ta-\aa} (v^2 + x^2 v_x^2) dxds
\end{split}
\end{equation}
for any $\oa>0$. To deal with $K_2$, one notes that
\begin{align*}
&(1+t)^{\varsigma}\int_{1/2}^1 \bar\rho^{\ga-\aa} \lt(x|r_x-1|+|r-x|\rt) |v| dx \\
\le &
(1+t)^\beta  \int_{1/2}^1 \bar\rho^{\ta-\aa} v^2 dx
+ (1+t)^{2\varsigma-\beta}\int_{1/2}^1 \bar\rho^{2\ga-\aa-\ta} \lt(x^2|r_x-1|^2+|r-x|^2\rt)dx \\
\le & C
(1+t)^\beta  \int_{1/2}^1 \bar\rho^{\ta} (  x^2 v_x^2 +v^2   )dx
+ (1+t)^{2\varsigma-\beta}\int_{1/2}^1 \bar\rho^{2\ga-\aa-\ta} \lt(x^2|r_x-1|^2+|r-x|^2\rt)dx.
\end{align*}
This implies that if $\ta+\aa\le \ga$,
\begin{align*}
&(1+t)^{\varsigma}\int_{1/2}^1 \bar\rho^{\ga-\aa} \lt(x|r_x-1|+|r-x|\rt) |v| dx \\
\le &
C(1+t)^\beta  \int  \bar\rho^{\ta} ( xv_x^2 + v^2 )dx
+ C(1+t)^{\beta-1}\int  \bar\rho^{\ga} \lt(x^2|r_x-1|^2+|r-x|^2\rt)dx;
\end{align*}
and if $\ta + \aa > \ga$,
\begin{align*}
&(1+t)^{\varsigma}\int_{1/2}^1 \bar\rho^{\ga-\aa} \lt(x|r_x-1|+|r-x|\rt) |v| dx \\
\le &
C(1+t)^\beta  \int  \bar\rho^{\ta} ( xv_x^2 + v^2 )dx + C \int  \bar\rho^{\ga-\aa} \lt(x^2|r_x-1|^2+|r-x|^2\rt)dx \\
& + C (1+t)^{\beta-1}\int  \bar\rho^{\ga} \lt(x^2|r_x-1|^2+|r-x|^2\rt)dx,
\end{align*}
due to
$$\frac{\aa}{\ga-\ta}(2\varsigma-\beta)\le \beta-1.$$
Then, one can obtain
\begin{equation}\label{estK2}\begin{split}
\int_0^t (1+s)^{\varsigma} | K_2 | ds \le  & C \int_0^t (1+s)^\beta  \int  \bar\rho^{\ta} ( xv_x^2 + v^2 )dx ds \\
&+ C \int_0^t \int  \bar\rho^{\ga-\aa} \lt(x^2|r_x-1|^2+|r-x|^2\rt)dx  ds \\
&  + C\int_0^t (1+s)^{\beta-1}\int  \bar\rho^{\ga} \lt(x^2|r_x-1|^2+|r-x|^2\rt)dxds.
\end{split}
\end{equation}
Note that for any $\oa>0$,
\begin{equation}\label{estK3}\begin{split}
\int_0^t (1+s)^{\varsigma} | K_3 | ds  \le & C \oa^{-1} \int_0^t (1+s)^{\varsigma}    \int  \bar\rho^{\ta} ( xv_x^2 + v^2 )dx ds \\
&+ \oa \int_0^t (1+s)^\varsigma \int  \bar\rho^{\ta -\aa} ( xv_x^2 + v^2 )dx ds .
\end{split}
\end{equation}
So, \ef{511-1} follows from \ef{lowest}, \ef{estlem51} and \ef{beauty}-\ef{estK3}.

{\em Step 3}.  In this step, we prove that
\begin{equation}\label{5-11-3}\begin{split}
(1+t)^{\frac{\beta+\varsigma}{2}  } \int\bar\rho^{\theta- \aa/2}(x^2v_{x}^2+v^2) (x,t) dx  \le C_\iota  \lt( \mathfrak{D}_1(0)+ \left\|  r_{0x}-1  \right\|_{L^\iy}^2  + \left\| v_x(\cdot, 0) \right\|_{L^\iy}^2\rt),
\end{split}\end{equation}
\begin{equation}\label{5-11-1}\begin{split}
(1+t)^\beta \int( v^2 + x^2 \bar\rho^{\theta} v_{x}^2) (x,t) dx  \le C_\iota  \lt( \mathfrak{D}_1(0)+ \left\|  r_{0x}-1  \right\|_{L^\iy}^2 + \left\| v_x(\cdot, 0) \right\|_{L^\iy}^2 \rt).
\end{split}\end{equation}
Notice that
\begin{align*}
&(1+t)^{\frac{\beta+\varsigma}{2}  } \int  \bar\rho^{\ta - \aa/2}  (x^2 v_x^2+v^2) (x, t)dx
\\
=& \int  \bar\rho^{\ta - \aa/2}  (x^2 v_x^2+v^2) (x, 0)dx
 +\frac{\beta+\varsigma}{2} \int_0^t (1+s)^{\frac{\beta+\varsigma}{2}-1} \int  \bar\rho^{\ta - \aa/2}  (x^2 v_x^2+v^2) dxds \\
& + 2 \int_0^t (1+s)^{\frac{\beta+\varsigma}{2}} \int  \bar\rho^{\ta - \aa/2}  (x^2 v_x v_{sx}+v v_s) dxds  \\
 \le & C \left\| v_x(\cdot, 0) \right\|_{L^\iy}^2 +  C\int_0^t (1+s)^{\beta} \int  \bar\rho^{\ta}  (x^2 v_x^2+v^2 + x^2 v_{sx}^2 + v_s^2) dxds \\
&+ C \int_0^t (1+s)^{ \varsigma} \int  \bar\rho^{\ta - \aa }  (x^2 v_x^2+v^2  ) dxds .
\end{align*}
Then, \ef{5-11-3} follows from \ef{estlem51}, \ef{estlem51'} and \ef{511-1}.
Similarly, one can show that
\begin{equation*}\label{}\begin{split}
(1+t)^\beta \int  \bar\rho^{\theta}(  x^2 v_{x}^2 + v^2 ) (x,t) dx  \le C_\iota  \lt( \mathfrak{D}_1(0)+ \left\|  r_{0x}-1  \right\|_{L^\iy}^2 + \left\| v_x(\cdot, 0) \right\|_{L^\iy}^2 \rt).
\end{split}\end{equation*}
This, together with \ef{phy}, \ef{hd1}, \ef{hd2} and $2(\ga-1)>\ga/2 \ge \ta$, gives \ef{5-11-1}.

\hfill$\Box$

Suppose that \ef{rx} and \ef{vx} hold.  Let $\ta\in (0,\ga/2]$, and
 $\alpha$ and $\beta$ be given  respectively in \ef{alpha} and \ef{beta}.  Then, for any $\iota\in (0, (2\ga-2-\ta)/8 ]$ and $l\in (0,1)$, one has

 \begin{lem}\label{lem3}  Suppose that \ef{rx} and \ef{vx} hold. Let $\ta\in (0,\ga/2]$, and
 $\alpha$, $\beta$  and $\varsigma$ be given  respectively in \ef{alpha}, \ef{beta} and \ef{varsa}. If $\iota\in (0, (2\ga-2-\ta)/8 ]$, then      for  $(x,t)\in I \times[0,T]$,
 \begin{align}
&(1+t)^{\frac{\ga-1+\aa-\theta}{\ga+\aa-\ta}\beta} x|r(x, t)-x|^2 +
(1+t)^{ \frac{3\beta+\varsigma}{4} - \frac{\beta-\varsigma}{4 \aa} \max\{0, \ 4\ta-4(\ga-1)-\aa\}   } x v^2(x,t)\notag\\
& +  x^3 | r_x(x, t)-1 |^2 + x^3 |v_x(x,t)|^2 \le  C_\iota  \lt( \mathfrak{D}_1(0)+ \left\|  r_{0x}-1  \right\|_{L^\iy}^2 + \left\| v_x(\cdot, 0) \right\|_{L^\iy}^2 \rt). \label{verify}
\end{align}
 \end{lem}

\noindent{\em Proof}. The proof consists of three steps.

{\em Step 1}.  In this step, we prove
\be\label{a3}
(1+t)^{\frac{\ga-1+\aa-\theta}{\ga+\aa-\ta}\beta} x|r(x, t)-x|^2
\le C_\iota  \lt( \mathfrak{D}_1(0)+ \left\|  r_{0x}-1  \right\|_{L^\iy}^2 \rt), \ee
\be\label{hha3}
(1+t)^{ \frac{3\beta+\varsigma}{4} - \frac{\beta-\varsigma}{4 \aa} \max\{0, \ 4\ta-4(\ga-1)-\aa\}   } x v^2(x,t)
\le C_\iota  \lt( \mathfrak{D}_1(0)+ \left\|  r_{0x}-1  \right\|_{L^\iy}^2 + \left\| v_x(\cdot, 0) \right\|_{L^\iy}^2 \rt). \ee

Clearly, \ef{a3}  follows from \ef{3.63} and
\be\label{a1} \begin{split}
& x(r(x, t)-x)^2= \int_0^x \lt(y(r(y, t)-y)^2\rt)_ydy\\
 \le & \int (r-x)^2(x, t) dx
+2\lt(\int   (r-x)^2dx\rt)^{1/2}\lt(\int x^2(r_x-1)^2dx\rt)^{1/2},
\ \ x\in I. \end{split} \ee
Similarly, for $ x\in I$, it holds that
\bee\label{} \begin{split}
 xv^2(x,t)= \int_0^x \lt(yv^2\rt)_ydy
 \le  \int v^2(x,t) dx
+2\lt(\int \bar\rho^{\aa/2-\ta}  v^2dx\rt)^{1/2}\lt(\int \bar\rho^{\ta-\aa/2} x^2 v_x^2 dx\rt)^{1/2}.\end{split} \eee
It follows from \ef{phy} and \ef{hd2} that
\bee\label{} \begin{split}
& \int \bar\rho^{\aa/2-\ta}  v^2dx \le C \int_0^{1/2}   v^2 dx
  + C \int_{1/2}^1 \bar\rho^{\aa/2-\ta + 2 (\ga-1)} ( v^2 + v_x^2 ) dx \\
  \le & C \int v^2 dx + C \lt(\int x^2 \bar\rho^{ \ta } v_x^2 dx \rt)^{\frac{2\aa-4\ta+4(
  \ga-1)}{\aa}} \lt(\int x^2 \bar\rho^{ \ta -\aa/2 }  v_x^2 dx \rt)^{\frac{4\ta-4(\ga-1)-\aa}{\aa}},
\end{split} \eee
if $\ta>\ga-1+\aa/4$, and if $\ta\le \ga-1+\aa/4$,
\bee\label{} \begin{split}
 \int \bar\rho^{\aa/2-\ta}  v^2dx
  \le   C \int v^2 dx + C  \int x^2 \bar\rho^{ \ta } v_x^2 dx .
\end{split} \eee
Then, \ef{hha3} follows from \ef{3.64}.

{\em Step  2}. In this step, we prove
\begin{equation}\label{lem3est}\begin{split}
 x^3 | r_x(x, t)-1 |^2 \le  C_\iota  \lt( \mathfrak{D}_1(0)+ \left\|  r_{0x}-1  \right\|_{L^\iy}^2 + \left\| v_x(\cdot, 0) \right\|_{L^\iy}^2 \rt).
\end{split}\end{equation}

Rewrite equation \eqref{nsp} as
\begin{equation}\label{ht7-1}\begin{split}
& \frac{\nu}{\ta} \lt\{\lt(\frac{r}{x}\rt)^{\frac{4\nu_1}{\nu}\ta} \lt[\bar\rho^\ta\lt(\frac{x^2}{r^2r_x}\rt)^\ta\rt]_x\rt\}_t \\
=& -\lt(\frac{r}{x}\rt)^{\frac{4\nu_1}{\nu}\ta} \lt\{\bar\rho\left( \frac{x}{r}\right)^2  v_t   +  \left[ \bar\rho^\ga \left(\frac{x^2}{r^2 r_x}\right)^\ga  \right]_x - \frac{x^4}{r^4}  \left(\bar{\rho}^\ga\right)_x\rt\}.
\end{split}
\end{equation}
Set
\be\label{Z}
Z(x, t)=\frac{\nu}{\ta}\lt(\frac{r}{x}\rt)^{\frac{4\nu_1\ta}{\nu}} \lt(\frac{x^2}{r^2r_x}\rt)^\ta-\frac{\nu}{\ta}- \bar\rho^{-\ta}(x)\int_x^1\bar\rho(y)\lt(\frac{r}{y}\rt)^{\frac{4\nu_1\ta}{\nu}-2}v(y,t)dy. \ee
Integrate  equation \ef{ht7-1} over $[x, 1]$ to get
\begin{equation}\label{ec1}\begin{split}
& \pl_t Z
 + \frac{\ta}{\nu}\bar\rho^{\ga-\ta}(x)\lt(\frac{x^2}{r^2r_x}\rt)^{\ga-\ta} Z+ \bar\rho^{\ga-\ta}(x) \lt[\lt(\frac{x^2}{r^2r_x}\rt)^{\ga-\ta} -\lt(\frac{r}{x}\rt)^{\frac{4\nu_1}{\nu}\ta-4} \rt] \\
  = & - \pl_t Y(x,t) -\sum_{i=1}^2 \mathfrak{L}_i(x, t) ,
\end{split}
\end{equation}
where $Y$ and $\mathfrak{L}_i$ satisfies the following estimates:
\begin{align*}
&Y(x, t)   =
 \frac{\nu}{\ta}\bar\rho^{-\ta}(x)  \int_x^1\bar\rho^{\ta}(y)\lt(\frac{y^2}{r^2r_y}\rt)^{\ta}
\lt(\lt(\frac{r}{y}\rt)^{\frac{4\nu_1\ta}{\nu}}\rt)_ydy , \\
&\mathfrak{L}_1(x, t)   =  \bar\rho^{-\ta}(x) \int_x^1 \bar\rho^{\ga}(y) \lt\{\lt[\lt(\frac{r}{y}\rt)^{\frac{4\nu_1}{\nu}\ta}\rt]_y  \lt(\frac{y^2}{r^2r_y}\rt)^{\ga} -\lt[\lt(\frac{r}{y}\rt)^{\frac{4\nu_1\ta}{\nu}-4} \rt]_y
\rt\}dy ,\\
& |\mathfrak{L}_2(x, t)| \le C \bar\rho^{\ga-2\ta}(x)\int_x^1 \bar\rho(y) |v(y, t)|dy
 +C\bar\rho^{-\ta}(x)\int_x^1  y^{-1} \bar\rho(y) v^2(y, t) dy.
\end{align*}
Rewrite $Z$ as
\be\label{defZ}
Z(x,t)=\frac{\nu}{\ta}\lt(r_x^{-\ta}-1\rt)+ Z_1(x,t) ,
\ee
where
\begin{align*}
Z_1(x,t)=   \frac{\nu}{\ta}\lt[\lt(\frac{r}{x}\rt)^{\frac{4\nu_1\ta}{\nu}-2\ta} -1\rt] r_x^{-\ta} - \bar\rho^{-\ta}(x)\int_x^1\bar\rho(y)\lt(\frac{r}{y}\rt)^{\frac{4\nu_1\ta}{\nu}-2}v(y,t)dy.
\end{align*}
Under the a priori assumptions \ef{rx}, we can see that the leading term of $Z$ is $\nu(1-r_x)$. Note that
$$
 \bar\rho^{\ga-\ta}(x) \lt[\lt(\frac{x^2}{r^2r_x}\rt)^{\ga-\ta} -\lt(\frac{r}{x}\rt)^{\frac{4\nu_1}{\nu}\ta-4} \rt]
 = \bar\rho^{\ga-\ta}(x) \lt( r_x^{\ta-\ga}-1\rt) +  \mathfrak{L} (x,t)  ,
$$
where
$$|\mathfrak{L} (x,t)| \le C x^{-1}\bar\rho^{\ga-\ta}(x) |r(x,t)-x|. $$
Then, \ef{rx} and the Taylor expansion imply that
 \be\label{ec2}\begin{split}
&\pl_t Z + \bar\rho^{\ga-\ta}(x)   a(x, t)Z
=  - \pl_t Y(x,t)  - \sum_{i=1}^3 \mathfrak{L}_i(x, t)  ,
\end{split}
 \ee
where
\begin{align*}
& \ga/(2\nu) \le a(x, t)\le  2 \ga /\nu,  \\
& |\mathfrak{L}_3(x, t)| \le C x^{-1}\bar\rho^{\ga-\ta}(x) |r(x,t)-x|
+ C \bar\rho^{\ga-2\ta}(x) \int_x^1 \bar\rho (y) |v(y,t)| dy,
\end{align*}
due to the smallness of $\ea_0$.
Integrate \ef{ec2} on $[0,t]$ to get
 \be\label{hnewec2}\begin{split}
Z(x,t) \le   Z(x,0)  + C\sup_{s\in [0,t]} |Y(x,s)|
 + \sum_{i=1}^3  \int_0^t  \exp\lt\{-\bar\rho^{\ga-\ta}(x) \int_s^t a(x,\tau) d\tau \rt\} |\mathfrak{L}_i(x, s)| ds  .
\end{split}
 \ee

Next, one needs to estimate the terms on the right-hand side of \ef{hnewec2}. It follows from $\aa > \ta$, $\ta \le \ga/2$ and \ef{phy} that
\be\label{hhYest}\begin{split}
& |Y(x, t) |  \le   C  \bar\rho^{-\ta}(x)  \int_x^1 y^{-2}\bar\rho^{\ta}(y)
\lt(y|r_y-1|+|r-y|\rt)dy \\
\le & C  \bar\rho^{-\ta}(x)  \lt(\int_x^1 y^{-4}\bar\rho^{ \ta+\aa}(y) dy\rt)^{1/2}
\lt(\int_0^1 \bar\rho^{ \ta- \aa}(y) \lt(y^2|r_y-1|^2+|r-y|^2\rt)dy \rt)^{1/2} \\
\le & C x^{-3/2}  \lt(\int_0^1 \bar\rho^{ \ta- \aa}(y) \lt(y^2|r_y-1|^2+|r-y|^2\rt)dy \rt)^{1/2},
\end{split}\ee
\be\label{hhL1est}\begin{split}
& |\mathfrak{L}_1(x, t) |  \le   C  \bar\rho^{-\ta}(x)  \int_x^1 y^{-2}\bar\rho^{\ga}(y)
\lt(y|r_y-1|+|r-y|\rt)dy \\
\le & C  \bar\rho^{-\ta}(x)  \lt(\int_x^1 y^{-4}\bar\rho^{2\ga-\ta+\aa}(y) dy\rt)^{1/2}
\lt(\int_0^1 \bar\rho^{ \ta- \aa}(y) \lt(y^2|r_y-1|^2+|r-y|^2\rt)dy \rt)^{1/2} \\
\le & C x^{-3/2} \bar\rho^{\ga-\ta}\lt(\int_0^1 \bar\rho^{ \ta- \aa}(y) \lt(y^2|r_y-1|^2+|r-y|^2\rt)dy \rt)^{1/2},
\end{split}\ee
\be\label{hhL2est}\begin{split}
 & |\mathfrak{L}_2(x, t)| \le C \bar\rho^{\ga-\ta}(x)\lt[ x^{-1/2} \lt(\int_0^1 y^2 \bar\rho v^2 dy \rt)^{1/2}+ x^{-1} \|(xv^2)(\cdot,t)\|_{L^\iy} \rt] ,
\end{split}\ee
\be\label{hhL3est}\begin{split}
 |\mathfrak{L}_3(x, t)| \le C \bar\rho^{\ga-\ta}(x) \lt[ x^{-1} |r(x,t)-x| + x^{-1/2} \lt(\int_0^1 y^2 \bar\rho v^2 dy \rt)^{1/2}  \rt] .
\end{split}\ee
Notice that
$$\bar\rho^{\ga-\ta}(x)\int_0^t \exp\lt\{-\bar\rho^{\ga-\ta}(x)\int_s^t a (x, \tau)d\tau\rt\}ds\le C$$
and
\be\label{hhZest}
|Z(x,t)|\le C  |r_x(x,t)-1|+ C x^{-1}|r(x,t)-x| + C x^{-1/2} \lt(\int_0^1 y^2 \bar\rho v^2 dy \rt)^{1/2}.
\ee
Then, \ef{lem3est} follows from \ef{hnewec2}, \ef{hhYest}-\ef{hhL2est}, \ef{estlem51}, \ef{a3} and \ef{hha3}.

{\em Step 3}. In this step, we prove
\begin{equation}\label{ptvx}\begin{split}
 x^3 |v_x(x, t)  |^2 \le  C_\iota  \lt( \mathfrak{D}_1(0)+ \left\|  r_{0x}-1  \right\|_{L^\iy}^2 + \left\| v_x(\cdot, 0) \right\|_{L^\iy}^2 \rt).
\end{split}\end{equation}
It follows from \ef{defZ}, \ef{rx} and \ef{vx} that
\bee\label{}
|v_x(x,t)|\le  C |\pl_t Z(x,t)| + C |\pl_t Z_1(x,t)|
 \eee
and
\bee\label{}\begin{split}
|\pl_t Z_1(x,t)| \le & C x^{-1}|v(x,t)| + C \ea_0 |v_x(x,t)| + C\bar\rho^{-\ta}(x)
\int_x^1 \bar\rho(y) (|v_t(y,t)|+|v(y,t)|) dy \\
\le & C x^{-1}|v(x,t)| + C \ea_0 |v_x(x,t)|  + C x^{-1/2} \lt(\int_0^1 y^2 \bar\rho (v^2+ v_t^2) dy \rt)^{1/2} .
\end{split}\eee
Due to the smallness of $\ea_0$, one has
\be\label{hh-vxest}
|v_x(x,t)|\le  C |\pl_t Z(x,t)| +  C x^{-1}|v(x,t)|   + C x^{-1/2} \lt(\int_0^1 y^2 \bar\rho (v^2+ v_t^2) dy \rt)^{1/2}.
 \ee
Note that
\begin{equation}\label{hh-dtY}\begin{split}
|\pl_t Y(x,t)|\le & C \bar\rho^{-\ta} \int_x^1 y^{-2}\bar\rho^\ta(y) (y|v_y|+|v|) dy \\
\le & C   \bar\rho^{-\ta} \lt( \int_x^1 y^{-4}\bar\rho^{2\ta-(\ga-1) }(y)  dy \rt)^{1/2} \lt( \int_0^1  \bar\rho^{\ga-1}(y) (y^2 v_y^2 + v^2) dy \rt)^{1/2} \\
\le &  C x^{-3/2}\lt( \int_0^1  \bar\rho^{\ga-1}(y) (y^2 v_y^2 + v^2) dy \rt)^{1/2} ,
  \end{split}\end{equation}
due to $\ta\le \ga/2 < 2(\ga-1)$. Then, \ef{ptvx} follows from \ef{hh-vxest}, \ef{ec2}, \ef{hhL1est}-\ef{hhZest}, \ef{a3}, \ef{hha3}, \ef{lem3est}, \ef{estlem51}, \ef{lowest} and \ef{3.64}.

 \hfill $\Box$

\subsection{Higher-order  estimates}\label{sect3.2}
For the higher-order estimates, we set
 \be \label{hhq} \mathfrak{Q}=\frac{x^2}{r^2r_x}-1.\ee
Then for any positive constant $k$,
\be\label{hh3.95}\begin{split}
\lt[ \bar\rho^{k}\left(\lt(\frac{x^2}{r^2r_x}\rt)^{k}-1\rt) \rt]_x
=\lt[ \bar\rho^{k} \lt( \lt( 1+  \mathfrak{Q}\rt)^{k}-1\rt) \rt]_x
= k \bar\rho^{k}  \lt( 1+  \mathfrak{Q}\rt)^{k-1}\mathfrak{Q}_x + \mathscr{P}_k(x,t);
\end{split}\ee
where $ \mathscr{P}_k$ satisfies the following estimate:
\be\label{hh3.96}
 |\mathscr{P}_k(x,t)| \le C x \bar\rho^{k-(\ga-1)}   \lt(|r_x-1|+ |r/x-1|\rt),
\ee
due to \ef{rhox}.
Using above notations and \ef{rhox}, one can rewrite \eqref{nsp} as
\begin{equation}\label{7-1}\begin{split}
 &\frac{\nu}{\ta}\lt[ \lt(\frac{r}{x}\rt)^{\frac{4\nu_1}{\nu}\ta} \lt(\ta \bar\rho^{\ta}  \lt( 1+  \mathfrak{Q}\rt)^{\ta-1}\mathfrak{Q}_x + \mathscr{P}_\ta\rt) \rt] _t +\lt(\frac{r}{x}\rt)^{\frac{4\nu_1}{\nu}\ta} \lt[ \ga \bar\rho^{\ga}  \lt( 1+  \mathfrak{Q}\rt)^{\ga-1}\mathfrak{Q}_x + \mathscr{P}_\ga \rt] \\
 & =\lt(\frac{r}{x}\rt)^{\frac{4\nu_1}{\nu}\ta} \lt[-\bar\rho\left( \frac{x}{r}\right)^2  v_t
+\lt(\frac{x^4}{r^4} -1\rt) \left(\bar{\rho}^\ga\right)_x\rt]
 +\frac{4\nu_1\ta}{\ga}\lt(\frac{r}{x}\rt)^{\frac{4\nu_1}{\nu}\ta-1} \phi\bar\rho^{\ta+1-\ga}v.
\end{split}
\end{equation}
The  principal part on the left-hand side of \ef{7-1} is
$${\nu} \bar\rho^{\ta} \mathfrak{Q}_{xt}+\ga \bar\rho^{\ga} \mathfrak{Q}_x = \nu(\bar{\rho}^\theta \mathfrak{Q}_x)_t + \gamma\bar{\rho}^{\gamma-\theta}(\bar{\rho}^\theta \mathfrak{Q}_x),$$
 which can be understood as a damped transport operator for $\bar\rho^\theta \mathfrak{Q}_x$ with a degenerate damping coefficient $\bar{\rho}^{\gamma-\theta}$. This is an interplay between the viscosity and the pressure. This structure leads to desirable estimates on the derivatives of $\mathfrak{Q}$, for example, the bound for
$$ \int \bar\rho^{2\ga -1} (x) \mathfrak{Q}_x^2(x, t)dx+\int_0^t\int  \bar\rho^{3\ga-1-\ta}  \mathfrak {Q}_x^2  dxds. $$
The bounds for the weighted norms of the derivatives of $\mathfrak{Q}$ in turn yield the bounds for those of the derivatives of $r$ and $v$, which are given by Lemma \ref{lem3.7} stated below.

In Lemma \ref{lem5}, one can use Lemma \ref{lem3.7} to bound
$$\int  \bar\rho^{2\ga-1}  \lt[ \lt( (r/x)_x \rt) ^2 + r_{xx}^2\rt](x,t)  dx$$
and to estimate the decay of
$$\int_0^l   \lt[ \lt( (r/x)_x \rt) ^2 + r_{xx}^2\rt](x,t)  dx , \ \  0<l<1. $$
 Based on this and other estimates, we can estimate the decay of
$$ \int \lt(\bar\rho v_t^2\rt) (x, t)dx  \ \ {\rm   and}    \ \  \int_0^l   \lt[ \lt( (v/x)_x \rt) ^2 + v_{xx}^2\rt](x,t)  dx ,  \ \  0<l<1.$$
 Putting those together yields the desired decay estimate for
 $$\lt\| (v_x, \ v/x, \ r_x-1, \ r/x-1)(\cdot, t) \rt\|_{H^1([0,l])}^2, \ \ 0<l<1 ,$$
 which, together with the supreme norm estimates for $r_x-1$ and $v_x$ away from the origin given in Lemma  \ref{lem3},
verifies the a priori assumptions \ef{rx} and \ef{vx} and closes the bootstrap argument.

We now carry out the strategy outlined above.

\begin{lem}\label{lem3.7}Suppose that \ef{rx} and \ef{vx} hold. Let $w(x)$ be a smooth function on $[0,1]$ satisfying $w(x)\ge 0$ and $w'(x)\le 0$ on $[0,1]$, and
$w(1)=0$. Then,
\begin{align}
&\int w \lt[4(r/x-1)^2 + (r_x-1)^2\rt] dx \le 2 \int w  \mathfrak{Q}^2  dx, \label{orirx}\\
&\int w \lt[  (v/x) ^2 + (v_x)^2\rt] dx \le C \int w  (\mathfrak{Q}_t)^2  dx,  \label{orivx}\\
&\int  w \lt[  \lt( (r/x)_x \rt)^2 + (r_{xx})^2\rt] dx \le C \int w ( \mathfrak{Q}_x )^2  dx,  \label{orirxx}\\
&\int  w \lt[  \lt((v/x)_x \rt)^2 + (v_{xx})^2\rt] dx \le C \int w \lt[ ( \mathfrak{Q}_{xt}  )^2 +    ( \mathfrak{Q}_x )^2  \rt]  dx.  \label{orivxx}
\end{align}
\end{lem}

\noindent{\em Proof}.  It follows from \ef{rx} that
$$\mathfrak{Q}= -2(r/x-1)-(r_x-1) + O(1) \ea_0 \lt(|r/x-1|+|r_x-1|\rt), $$
which, together with the integration by parts, $r(0,t)=0$ and $w'(x)\le 0$, gives
\begin{equation}\label{}\begin{split}
\int w  \mathfrak{Q}^2  dx  \ge & \int w \lt[4(r/x-1)^2 + (r_x-1)^2\rt] dx -  2 \int (  x^{-1} w )_x  (r-x)^2  dx \\
 &- C \ea_0
\int w \lt[ (r/x-1)^2 + (r_x-1)^2\rt] dx \\
\ge & \frac{1}{2} \int w \lt[4(r/x-1)^2 + (r_x-1)^2\rt] dx,
\end{split}
\end{equation}
due to $(  x^{-1} w)_x \le 0$ and the smallness of $\ea_0$. This proves \ef{orirx}. Similarly, \ef{orivx} follows from
\begin{align}\label{hhQt}
\mathfrak{Q}_t = - (1+ \mathfrak{Q}) \lt(2v/r+ v_x/r_x\rt).
\end{align}
For \ef{orirxx}, note that
\begin{align}\label{hhQx}
\mathfrak{Q}_x = - (1+ \mathfrak{Q}) \lt[2 (r/x)(r/x)_x+ r_{xx}/r_x\rt].
\end{align}
Thus,
\begin{align}\label{}
\int w (\mathfrak{Q}_x)^2 \ge & 2c \int  w   \lt[2(r/x)_x+ r_{xx} \rt]^2 dx - C\ea_0 \int  w   \lt[ \lt((r/x)_x\rt)^2+ (r_{xx})^2 \rt]  dx \notag\\
 \ge & c\int  w   \lt[ \lt((r/x)_x\rt)^2+ (r_{xx})^2 \rt]  dx, \notag
\end{align}
due to the smallness of $\ea_0$, and
\begin{align}\label{}
4\int  w    (r/x)_x  r_{xx}  dx
= &8 \int  w \lt((r/x)_x\rt)^2 dx
-2 \int  (   x w )_x  \lt((r/x)_x\rt)^2 dx \notag \\
= & 6 \int  w \lt((r/x)_x\rt)^2 dx
-2 \int  xw_x   \lt((r/x)_x\rt)^2 dx \ge 0. \notag
\end{align}
Similarly, one can prove \ef{orivxx}.

\hfill $\Box$

\begin{lem}\label{lem5} Suppose that \ef{rx} and \ef{vx} hold.  Let $\ta\in (0,\ga/2]$, and
 $\alpha$ and $\beta$ be given  respectively in \ef{alpha} and \ef{beta}.  If $\iota\in (0, (2\ga-2-\ta)/8 ]$ and $l\in (0,1)$, then for $t\in [0,T]$,
\begin{align}
&  \int  \bar\rho^{2\ga-1}  \lt[ \lt( (r/x)_x \rt) ^2 + r_{xx}^2\rt](x,t)  dx
   \le  C \mathscr{E}(0),  \label{rxx-est} \\
& (1+t)^{\beta-1} \int_0^l   \lt[ \lt( (r/x)_x \rt) ^2 + r_{xx}^2\rt](x,t)  dx    \le  C_{\iota,l} \mathscr{E}(0), \label{decayrxx} \\
& \int_0^t (1+s)^{\beta-1}   \int \bar\rho^\ta \lt(v_x^2+(v/x)^2\rt) dxds\le C_\iota \mathscr{E}(0) .  \label{hh155}
\end{align}
\end{lem}

\noindent{\em Proof}. The proof consists three steps.

{\em Step 1}.  In this step, we prove
\be\label{hh12172}\begin{split}
\int \bar\rho^{2\ga -1} (x) \mathfrak{Q}_x^2(x, t)dx+\int_0^t\int  \bar\rho^{3\ga-1-\ta}  \mathfrak {Q}_x^2  dxds
 \le C \mathscr{E}(0).
\end{split}\ee
This, together with \ef{orirxx}, gives \ef{rxx-est}.

It follows from \ef{7-1}, \ef{rhox} and the Cauchy inequality that
\begin{equation}\label{ht3104}\begin{split}
&\int \bar\rho^{2\ga-1}   | \mathfrak{Q}_x(x,t)|^2 dx
+\int_0^t \int \bar\rho^{3\ga-1-\ta}   | \mathfrak{Q}_x |^2 dx
  \\
\le & C  \mathscr{E}(0) + C \int \bar\rho(x) \lt[x^2(r_x-1)^2 + (r-x)^2\rt](x,t) dx +C \int_0^t \int \bar\rho^{\ta -(\ga-1)} v^2 dx ds \\
&+C \int_0^t \int \bar\rho^{\ga+1-\ta}  \lt[x^2(r_x-1)^2 + (r-x)^2 + v_s^2\rt] dx ds .
\end{split}
\end{equation}
Notice that
\begin{equation*}\label{}\begin{split}
\int \bar\rho^{\ta-(\ga-1)} v^2 dx \le & C \int_0^{1/2} x^2 (v^2+v_x^2) + C \int_{1/2}^1  \bar\rho^{\ta+(\ga-1)} (v^2+ v_x^2) dx \\
 \le & C \int \bar\rho^\ta (v^2+ x^2 v_x^2) dx ,
\end{split}
\end{equation*}
where \ef{hd1}, \ef{hd2} and \ef{phy} have been used. Then,\ef{hh12172} follows from \ef{lowest} and $\ta \le 1<\ga$.

{\em Step 2}. In this step, we prove \ef{decayrxx}.
Let $\psi$ be a smooth cut-off function defined on $[0,1]$ satisfying
\begin{equation}\label{psi}\begin{split}
 \psi=1 \ \ {\rm on } \ \ [0,1-l], \ \ \psi=0 \ \ {\rm on } \ \  [1-l/2, 1] \ \ {\rm and} \  \ -  8/ l  \le  \psi'(x)\le 0\  \ {\rm on } \ \  [0, 1]
\end{split}
\end{equation}
for any fixed constant $l\in (0,1)$.
As shown in Step 1, one can obatin
\begin{align}
&  \frac{\nu}{\ta}\frac{d}{dt}\int \psi\bar\rho^{2\ga-1-2\ta}\lt[ \lt(\frac{r}{x}\rt)^{\frac{4\nu_1}{\nu}\ta} \lt(\ta \bar\rho^{\ta}  \lt( 1+  \mathfrak{Q}\rt)^{\ta-1}\mathfrak{Q}_x + \mathscr{P}_\ta\rt) \rt]^2 dx
+  \int \psi \bar\rho^{3\ga-1-\ta}   | \mathfrak{Q}_x |^2 dx
 \notag \\
& \le    C  \int \bar\rho^{\ta } ( v^2 + x^2 v_x^2 + v_s^2 )dx
 +C   \int \bar\rho^{\ga}  \lt[x^2(r_x-1)^2 + (r-x)^2 \rt] dx ,
\end{align}
which, together with \ef{estlem51}, \ef{3.63} and \ef{hh12172}, implies that
\be\label{haha.11}\begin{split}
(1+t)^{\beta-1}\int \psi \bar\rho^{2\ga -1} (x) \mathfrak{Q}_x^2(x, t)dx+\int_0^t (1+s)^{\beta-1}\int  \psi \bar\rho^{3\ga-1-\ta} \mathfrak {Q}_x^2dxds
 \le C_{\iota,l} \mathscr{E}(0).
\end{split}\ee
(If $\beta>2$, we first show the time decay with rate $-1$. With this, we can then show the time decay with rate $1-\beta$. Indeed, this technique has been used in Step 3 of the proof of Lemma \ref{lem3.4}.)

Squaring \ef{7-1}, multiplying the resulting equation by $\psi \bar\rho^{3\ga-1-3\ta}$, and  integrating the product over $[0, 1] $, one obatins that
\begin{align}
&\int_0^t (1+s)^{\beta-1} \int \psi \bar\rho^{3\ga-1-\ta} \mathfrak{Q}_{xs}^2 dx ds
\le   C \int_0^t  (1+s)^{\beta-1} \int \psi \bar\rho^{3\ga-1-\ta} \mathfrak{Q}_{x}^2 dxds   \notag\\
& +C \int_0^t  (1+s)^{\beta-1} \int \psi \bar\rho^{\ga+1-\ta}  \lt[x^2(r_x-1)^2 + (r-x)^2 + v_s^2 + x^2 v_x^2 +v^2 \rt] dx ds, \label{hr108}
\end{align}
where \ef{vx} has been used.  This, together with \ef{haha.11},  \ef{estlem51}, \ef{3.63} and $\ta\le 1<\ga$, gives
\begin{align}\label{hn3.98}
\int_0^t (1+s)^{\beta-1} \int \psi \bar\rho^{3\ga-1-\ta} \mathfrak{Q}_{xs}^2 dx ds \le  C_{\iota,l} \mathscr{E}(0).
\end{align}
So, it follows from \ef{orirxx} \ef{orivxx},  \ef{haha.11} and \ef{hn3.98} that
\begin{align}
& (1+t)^{\beta-1} \int \psi \bar\rho^{2\ga-1}  \lt[ \lt( (r/x)_x \rt) ^2 + r_{xx}^2\rt](x,t)  dx \notag\\
& +
\int_0^t (1+s)^{\beta-1} \int \psi \bar\rho^{3\ga-1-\ta}  \lt[ \lt( (v/x)_x \rt) ^2 + v_{xx}^2\rt]  dx ds \le  C_{\iota,l} \mathscr{E}(0). \label{ho}
\end{align}

{\em Step 3}. In this step, we prove \ef{hh155}. Choose $l=1/2$ in \ef{ho} to get
\be\label{xuyao}\int_0^t (1+s)^{\beta-1} \int_0^{1/2}  \lt[ \lt( (v/x)_x \rt) ^2 + v_{xx}^2\rt]  dx ds \le  C_{\iota} \mathscr{E}(0).
\ee
Note that
\begin{align}
\int_0^{1/2} \lt(v_x^2+(v/x)^2\rt) dx \le & C \int_0^{1/2} x^2 \lt[ v_x^2 + v_{xx}^2 +(v/x)^2 +\lt((v/x)_x\rt)^2 \rt] \notag\\
\le &  C \int_0^{1/2}  \lt[  x^2 v_x^2 + v^2 + v_{xx}^2  +\lt((v/x)_x\rt)^2 \rt], \label{smiles}
\end{align}
due to  \ef{hd1}. Then,  \ef{hh155} follows from \ef{estlem51} and \ef{xuyao}.

\hfill $\Box$

\begin{lem}\label{lem5-new} Suppose that \ef{rx} and \ef{vx} hold.  Let $\ta\in (0,\ga/2]$, and
 $\alpha$ and $\beta$ be given  respectively in \ef{alpha} and \ef{beta}.  If $\iota\in (0, (2\ga-2-\ta)/8 ]$ and $l\in (0,1)$, then for $t\in[0,T]$,
\begin{equation}\label{kd4}\begin{split}
(1+t)^{\beta-1} \int \lt(\bar\rho v_t^2\rt) (x, t)dx
+ \int_0^t (1+s)^{\beta-1} \int \bar\rho^{\ta} \left( v_{xs}^2
  + ( v_s/x)^2  \right) dxds
  \le  C_{\iota} \mathscr{E} (0),
\end{split}
\end{equation}
\be\label{kd7}
(1+t)^{\beta-1} \int_0^l   \lt[ \lt( (v/x)_x \rt) ^2 + v_{xx}^2\rt](x,t)  dx    \le  C_{\iota,l} \mathscr{E}(0).
 \ee
\end{lem}

\noindent{\em Proof}.  The proof consists of two steps. We prove \ef{kd4} and \ef{kd7} in Steps 1 and 2, respectively.

{\em Step 1}. In this step, we prove that
\begin{align}
 (1+t)^{\beta-1} \int  \lt( \psi \bar\rho v_t^2\rt) (x, t)dx
 + \int_0^t (1+s)^{\beta-1} \int \psi \bar\rho^{\ta} \left( v_{xs}^2
  + ( v_s/x)^2  \right) dxds
  \le  C_{\iota,l} \mathscr{E} (0), \label{shi}
\end{align}
where  $\psi$ is a smooth cut-off function defined on $[0,1]$ satisfying \ef{psi}. Then, \ef{kd4} follows from  \ef{3.63} and  \ef{shi}.

Differentiating \eqref{nsp} with respect to $t$ yields that
\begin{align}\label{nsptime}
 &\bar\rho\left( \frac{x}{r}\right)^2  v_{tt}  -2\bar\rho\left( \frac{x}{r}\right)^3 \frac{v}{x}  v_{t}  -\ga \left[ \rho^\ga \left(2\frac{v}{r}+\frac{v_x}{r_x}\right)   \right]_{x} +4 \left( \frac{x}{r}\right)^5 \frac{v}{x}\left(\bar{\rho}^\ga\right)_x   \notag \\
 = & \nu \lt[\rho^{\ta} \left(\frac{v_{xt}}{r_x}+2\frac{v_t}{r}\right)\rt]_x -  \nu \lt[\rho^{\ta} \left( \frac{v_x^2 }{r_x^2} +2 \frac{v^2}{r^2}  \right)\rt]_x-\ta\nu\lt[ \rho^{\ta}\left(\frac{v_{x}}{r_x}+2\frac{v}{r}\right)^2\rt]_x
\notag\\
&+4\nu_1\ta\lt[\rho^{\ta}\left(\frac{v_{x}}{r_x}+2\frac{v}{r}\right)\rt]_x\frac{v}{r} -4\nu_1(\rho^{\ta})_x\lt( \frac{v_t}{r} - \frac{v^2}{r^2} \rt),
\end{align}
where
$\rho = \rho(r(x,t),t) = \bar\rho r_x^{-1}(r/x)^{-2}$.
Multiplying   equation \eqref{nsptime} by $ \psi v_t$,  integrating the product with respect to  the  spatial variable, and noting that $v_t(0, t)=v(0, t)=0$,   one has
 \begin{equation}\label{levt'}\begin{split}
 \frac{d}{dt}\int \frac{1}{2}\bar{\rho} \psi \left(\frac{x}{r}\right)^2  v_t^2 dx
+ \nu  \int \rho^{\ta} \left(\frac{v_{xt}}{r_x}+2\frac{v_t}{r}\right)   \left(\psi v_t\right)_x dx
=\sum_{k=1}^5 J_k,
\end{split}
\end{equation}
where
\begin{equation*}\label{}\begin{split}
J_1=&\int \frac{v}{x} \bar{\rho}\psi\left(\frac{x}{r}\right)^3   v_{t}^2dx + 4 \int \left(\frac{x}{r}\right)^5v  \phi \bar{\rho}\psi v_t dx \le C
 \int \bar\rho^\ta  \left( v^2 + v_t^2\right)  dx,\\
J_2=&-\ga \int   \rho^\ga \left(2\frac{v}{r}+\frac{v_x}{r_x}\right)   \left(\psi v_t\right)_x dx,   \\
J_3= & \nu   \int  \rho^{\ta}\left[  \frac{v_x^2 }{r_x^2}  +2 \frac{v^2}{r^2}  + \ta \left(  \frac{v_x }{r_x}  +2 \frac{v}{r}    \right)^2  \right] \left(\psi v_t\right)_x dx , \\
J_4=& -4\nu_1\ta\int \rho^{\ta}\left(  \frac{v_x }{r_x}  +2 \frac{v}{r}    \right) \left(\psi v_t \frac{v}{r}\right)_x dx, \\
J_5= &4\nu_1\int \rho^{\ta}\left[  \lt( \frac{v_t}{r} - \frac{v^2}{r^2} \rt)_{x} \psi v_t  + \lt( \frac{v_t}{r} - \frac{v^2}{r^2} \rt)(\psi v_t)_x    \right]dx.
\end{split}
\end{equation*}

In view of \ef{hh3.95}, $\psi'(x) \le 0 $  and $(\bar\rho^\ta)_x \le 0$, the second term on the left-hand side of \eqref{levt'} can be estimated as follows:
 \begin{equation}\label{kd}\begin{split}
 &\int  \rho^{\ta}\left(\frac{v_{xt}}{r_x}+2\frac{v_t}{r}\right)  \left(\psi v_t\right)_x  dx \\
 =& \int \rho^{\ta}\psi\frac{v_{xt}^2}{r_x} dx
 +2 \int \rho^{\ta} \frac{\psi}{r} v_t v_{tx}dx + \int \rho^{\ta}\psi' \left(\frac{v_{xt}}{r_x}+2\frac{v_t}{r}\right) v_t  dx \\
 =& \int \rho^{\ta}\psi \lt(\frac{v_{xt}^2}{r_x} +r_x \frac{v_t^2}{r^2}  \rt)dx   - H( t)  + \int \rho^{\ta}\psi' \left(\frac{v_{xt}}{r_x}+2\frac{v_t}{r}\right) v_t  dx,
\end{split}
\end{equation}
where
$$H(t)=\int\frac{ \psi'(x)\rho^{\ta}}{r}v_t^2dx+\int \frac{\psi}{r} \lt[ \lt( \bar\rho^{\ta}\rt)_x+\mathcal{W}_x\rt]v_t^2  dx \le  \int \frac{\psi}{r} \mathcal{W}_x v_t^2 dx, $$
and
$$\mathcal{W}_x =\lt[ \bar\rho^{\ta}\left(\lt(\frac{x^2}{r^2r_x}\rt)^{\ta}-1\rt) \rt]_x
= \ta \bar\rho^{\ta}  \lt( 1+  \mathfrak{Q}\rt)^{\ta-1}\mathfrak{Q}_x + \mathscr{P}_\ta(x,t).
$$
It follows from \ef{hh12172}, \ef{lem2est'}, and the H\"{o}lder and Cauchy-Schwarz  inequalities that  for any $\oa>0$,
\begin{equation*}\begin{split}
 & H(t) \le \frac{1}{4}\int\frac{ \psi\rho^{\ta}r_x}{r^2}v_t^2dx+ C\int \psi \mathcal{W}_x^2 v_t^2dx\\
=&\frac{1}{4}\int\frac{ \psi\rho^{\ta}r_x}{r^2}v_t^2dx + 2C\int_0^{1-l/2} \psi \mathcal{W}_x^2 \lt(\int_0^x v_tv_{ t y} dy \rt)dx \\
= & \frac{1}{4}\int\frac{ \psi\rho^{\ta}r_x}{r^2}v_t^2dx + C \lt(\int  \psi \mathcal{W}_x^2 dx \rt) \lt( \int_0^{1-l/2} y^2v_{ty}^2dy\rt)^{1/2} \lt( \int_0^{1-l/2}\frac{v_t^2}{y^2}dy\rt)^{1/2}\\
\le & \frac{1}{4}\int\frac{ \psi\rho^{\ta}r_x}{r^2}v_t^2dx+ C_l \oa \int_0^{1-l}\psi\frac{v_t^2}{x^2}dx + C_l\oa\int_{1-l}^{1-l/2}
v_t^2 dx + C_l\oa^{-1}\int_0^{1-l/2}x^2v_{xt}^2dx\\
\le & \frac{1}{2}\int\frac{ \psi\rho^{\ta}r_x}{r^2}v_t^2dx + C_l\int \bar\rho^\ta ( x^2v_{xt}^2+  v_t^2) dx .
 \end{split}\end{equation*}
This, together with \ef{kd}, implies
\begin{equation}\label{kd1}\begin{split}
 \int  \rho^{\ta}\left(\frac{v_{xt}}{r_x}+2\frac{v_t}{r}\right)  \left(\psi v_t\right)_x dx
 \ge  c\int \psi\lt({v_{xt}^2} + \frac{ v_t^2}{x^2}\rt) dx
-C_l \int  \bar\rho^\ta \left( {x^2v_{xt}^2}  +   v_t^2\right)  dx.
\end{split}
\end{equation}
Thus, by \ef{levt'}, one has
 \begin{equation*}\label{kd3}\begin{split}
 \frac{d}{dt}\int \bar{\rho}\psi \left(\frac{x}{r}\right)^2  v_t^2 dx
+  \int \psi   \left(  {v_{xt}^2}
  + \frac{v_t^2}{x^2} \right)dx
\le   C_l
 \int  \bar\rho^\ta \left( {x^2v_{xt}^2}  +   v_t^2+v^2\right)  dx +  \sum_{k=2}^6 J_k.
\end{split}
\end{equation*}
It therefore follows easily from the Cauchy-Schwarz inequality that
 \begin{equation}\label{wawa2}\begin{split}
 \frac{d}{dt}\int \frac{1}{2}\bar{\rho}\psi \left(\frac{x}{r}\right)^2  v_t^2 dx
+  \int \psi \bar\rho^{\ta}\left(  {v_{xt}^2}
  + \frac{v_t^2}{x^2} \right)dx  \\
\le    C_l
 \int  \bar\rho^\ta \left( {x^2v_{xt}^2}  +   v_t^2+( v/x)^2+ v_x^2 \right)  dx .
\end{split}
\end{equation}
So, \ef{shi} follows from \ef{wawa2}, \ef{3.63} and \ef{hh155}.

{\em Step 2}. In a similar way to deriving \ef{hn3.98}, one can use \ef{haha.11}, \ef{shi}, \ef{estlem51} and \ef{3.63}
to obtain
$$(1+t)^{\beta-1}\int \psi (x) \bar\rho^{2\ga -1} (x) \mathfrak{Q}_{tx}^2(x, t)dx
\le   C_{\iota,l} \mathscr{E} (0).$$
This, together with \ef{orivxx} and \ef{haha.11}, gives \ef{kd7}.

\hfill  $\Box$

\begin{lem}\label{lem1521}
Suppose that \ef{rx} and \ef{vx} hold.  Let $\ta\in (0,\ga/2]$, and
 $\alpha$ and $\beta$ be given  respectively in \ef{alpha} and \ef{beta}.  If   $\iota\in (0, (2\ga-2-\ta)/8 ]$ and $l\in (0,1)$, then for $t\in [0,T]$,
\begin{align}
(1+t)^{\beta-1} \lt\| (v_x, \ v/x, \ r_x-1, \ r/x-1)(\cdot, t) \rt\|_{H^1([0,l])}^2  \le C_{\iota,l} \mathscr{E} (0) . \label{H1norm}
\end{align}
\end{lem}

\noindent{\em Proof}. In a similar way as for \ef{smiles}, one can use \ef{estlem51}, \ef{3.63}, \ef{decayrxx} and \ef{kd7} to obtain
\begin{align}
(1+t)^{\beta-1}\int_0^{1/2} \lt[v_x^2+(v/x)^2 +(r_x-1)^2+(r/x-1)^2 \rt](x,t) dx \le C_{\iota} \mathscr{E} (0) , \label{}
\end{align}
which implies, by using \ef{estlem51} and \ef{3.63} again, that for any $l\in (0,1)$,
\begin{align}
(1+t)^{\beta-1}\int_0^{l} \lt[v_x^2+(v/x)^2 +(r_x-1)^2+(r/x-1)^2 \rt](x,t) dx \le C_{\iota,l} \mathscr{E} (0) . \label{}
\end{align}
This, together with \ef{decayrxx} and \ef{kd7}, gives  \ef{H1norm}.

\hfill $\Box$

\subsection{Global existence of strong solutions}
In this subsection, we prove Theorem \ref{hhmainthm1} by first verifying the a priori assumptions \ef{rx} and \ef{vx}. Choose $\iota=  (2\ga-2-\ta)/8$ and $l=1/2$ in \ef{H1norm} to get
$$\lt\| (v_x, \ v/x, \ r_x-1, \ r/x-1)(\cdot, t) \rt\|_{H^1([0,1/2])}^2  \le C  \mathscr{E} (0), \ \  t\in [0,T],$$
which implies, using $\|\cdot\|_{L^\iy} \le \|\cdot\|_{H^1}$, that
$$\lt\| (v_x, \ v/x, \ r_x-1, \ r/x-1)(\cdot, t) \rt\|_{L^\iy([0,1/2])}^2  \le C  \mathscr{E} (0), \ \  t\in [0, T].$$
It follows from choosing  $\iota=  (2\ga-2-\ta)/8$ in \ef{verify} that
$$\lt\| (v_x, \ v/x, \ r_x-1, \ r/x-1)(\cdot, t) \rt\|_{L^\iy([1/2,1])}^2  \le C  \mathscr{E} (0), \ \  t\in [0, T].$$
This verifies the a priori assumption for small $\mathscr{E} (0)$.

The local existence of strong solutions for the problem \ef{419} can be obtained easily by combining the approximation techniques
in \cite{jangnsp} and the a priori estimates obtained in Sections \ref{sect3.1} and \ref{sect3.2}, at least for small $\mathscr{E} (0)$. Indeed, the a priori estimates obtained in Sections \ref{sect3.1} and \ref{sect3.2} are also sufficient for the local existence theory, for small $\mathscr{E} (0)$. Therefore, a standard continuation argument proves Theorem \ref{hhmainthm1}, with the estimates obtain in Sections \ref{sect3.1} and \ref{sect3.2}.

\begin{rmk} The local existence theory in \cite{duan} which is for the case of  $\ta=1$, does not apply to our case. \end{rmk}

\section{Nonlinear Asymptotic Stability}
This section is devoted to proving  Theorem \ref{mainthm2}.

First, it  follows from the fact that
$\|g\|_{L^\iy}^2 \le \| g\|_{L^2}\| g_x\|_{L^2}$ for any function $g$, \ef{3.63}, \ef{3.64}, \ef{H1norm} and \ef{verify} that
\be\label{decayofr-est}
 (1+t)^{\min\lt\{\frac{3\ga-2+2(\aa-\theta)}{2(\ga+\aa-\ta)}\beta
 -\frac{1}{2}   , \ \frac{\ga-1+\aa-\theta}{\ga+\aa-\ta}\beta  \rt\} }   |r(x, t)-x|^2 \le C_\iota \mathscr{E}(0)  , \ \ x\in I,
\ee
\be\label{decayofv-est}
\lt[(1+t)^{\frac{\beta}{2}}+ x(1+t)^{ \frac{3\beta+\varsigma}{4} - \frac{\beta-\varsigma}{4 \aa} \max\{0, \ 4\ta-4(\ga-1)-\aa\}   }   \rt]|u(r(x, t), t)|^2 \le C_\iota \mathscr{E}(0)  , \ \ x\in I. \ee
Since for small $\iota$,
$$\frac{3\ga-2+2(\aa-\theta)}{2(\ga+\aa-\ta)}\beta-\frac{1}{2} \ge \frac{\ga-1+\aa-\theta}{\ga+\aa-\ta}\beta$$
then we have \ef{decayofr} and \ef{decayofv}.

Due to the fact that $\|\cdot\|_{L^\iy} \le \|\cdot\|_{H^1}$ and \ef{H1norm}, we have for any fixed $l\in(0,1)$,
\begin{align}
& (1+t)^{\beta-1}   \lt(|v_x(x,t)|^2 + |x^{-1}v(x,t)|^2  + |r_x(x,t)-1|^2 + |x^{-1}r(x,t)-1|^2 \rt) \notag\\
& \le   C_{\iota, l}  \mathscr{E} (0), \ \ \ \  x\in [0,l], \label{sec43}
\end{align}
which implies \ef{decayofrx}.
To prove \ef{decayforrho}, we notice that for any $b\in [0, 2-\ga]$,
\begin{align}
&x^3\bar\rho^{-b}(x)|\rho(r(x,t),t)-\bar\rho(x)|^2 = x^3  \bar\rho^{2-b}(x) \mathfrak{Q}^2(x,t) \notag\\
 =& \int_0^x (y^3 \bar\rho^{2-b} \mathfrak{Q}^2)_y dy
 \le   3 \int_0^x y^2 \bar\rho^{2-b} \mathfrak{Q}^2 dy
 + 2 \int_0^x y^3 \bar\rho^{2-b} \mathfrak{Q}\mathfrak{Q}_y dy \notag\\
 \le &  C \int_0^1 y^2 \bar\rho^\ga (y^2|r_y-1|^2 + |r-y|^2)dy \notag\\
& + C \lt( \int \bar\rho^{2\ga-1}  \mathfrak{Q}_y^2 dy\rt)^{1/2} \lt( \int y^2 \bar\rho^{5-2\ga-2b}  \mathfrak{Q}^2 dy\rt)^{1/2}, \notag
\end{align}
where the first inequality is due to \ef{rhox}; and if $\ga>(5-2b)/3$
\begin{align}
  \int y^2 \bar\rho^{5-2b-2\ga}  \mathfrak{Q}^2 dy
  \le &  \lt(\int  y^2 \bar\rho^\ga (y^2|r_y-1|^2 + |r-y|^2)dy\rt)^{\frac{5-2\ga-\ta-2b}{\ga-\ta}} \notag \\
  & \times\lt(\int  y^2 \bar\rho^\ta (y^2|r_y-1|^2 + |r-y|^2)dy\rt)^{\frac{3\ga-5+2b}{\ga-\ta}}. \notag
\end{align}
Then, it follows from \ef{estlem51} and \ef{hh12172} that for any $b\in [0, 2-\ga]$,
\begin{align}
&x^3 \bar\rho^{-b}(x)(1+t)^{\frac{\beta}{2}-\frac{\max\{0, \ 3\ga-5+2b\}}{2(\ga-\ta)}} |\rho(r(x,t),t)-\bar\rho(x)|^2  \le  C  \mathscr{E} (0), \ \  x\in [0,1]. \label{sec41}
\end{align}
Due to \ef{sec43} and $b<1$, we have
\begin{align}
(1+t)^{\beta-1} \bar\rho^{-b}(x) |\rho(r(x,t),t)-\bar\rho(x)|^2  \le  C_{\iota, l}  \mathscr{E} (0), \ \  x\in [0,l]. \label{sec42}
\end{align}
So, \ef{decayforrho} is a consequence of \ef{sec41} and \ef{sec42}.

It remains to proving \ef{deforvx}. We use \ef{hh-vxest}, \ef{ec2}, \ef{hhL2est}-\ef{hhZest},  \ef{hh-dtY} and
\bee\label{}\begin{split}
  |\mathfrak{L}_1(x, t) |  \le   &
   C x^{-3/2} \lt(\int_0^1 \bar\rho^{ 2\ga - 2\ta +(\ga-1)}(y) \lt(y^2|r_y-1|^2+|r-y|^2\rt)dy \rt)^{1/2} \\
   \le &  C x^{-3/2} \lt(\int_0^1 \bar\rho^{ \ga }(y) \lt(y^2|r_y-1|^2+|r-y|^2\rt)dy \rt)^{1/2}
\end{split}\eee
 to get
\be\label{sec4.4}\begin{split}
& x^3|v_x(x,t)|^2 \le  C x^3 \bar\rho^{2\ga-2\ta}(x)| r_x(x,t)-1|^2
 +  C \int y^2 \bar\rho^{\ga-1} v_y^2(y,t) dy \\
& \quad +  C_{\iota}   \lt[(1+t)^{-\frac{\ga-1+\aa-\theta}{\ga+\aa-\ta}\beta} +  (1+t)^{- \frac{3\beta+\varsigma}{4} + \frac{\beta-\varsigma}{4 \aa} \max\{0, \ 4\ta-4(\ga-1)-\aa\}   } \rt] \mathscr{E} (0).
\end{split}\ee
Here \ef{3.63}, \ef{3.64} and \ef{estlem51}  have been used. In view of \ef{3.64}, we see that if $\ga-1\ge \ta$,
$$\int y^2 \bar\rho^{\ga-1} v_y^2(y,t) dy \le C_{\iota}(1+t)^{-\beta}\mathscr{E} (0);$$
and if $\ga-1<\ta$,
\bee\label{}\begin{split}
\int y^2 \bar\rho^{\ga-1} v_y^2(y,t) dy \le  &C \lt(\int y^2 \bar\rho^{\ta} v_y^2(y,t) dy\rt)^{\frac{\aa-2\ta+2(\ga-1)}{\aa}}\lt(\int y^2 \bar\rho^{\ta-\aa/2} v_y^2(y,t) dy\rt)^{\frac{ 2\ta-2(\ga-1)}{\aa}} \\
 \le & C_{\iota}(1+t)^{-\beta+ \frac{\beta-\varsigma}{2\aa}(2\ta-2(\ga-1))}\mathscr{E} (0).
\end{split}\eee
That means
\be\label{sec4.5}\begin{split}
\int y^2 \bar\rho^{\ga-1} v_y^2(y,t) dy
 \le & C_{\iota}(1+t)^{-\beta+ \frac{\beta-\varsigma}{2\aa}\max\{0, \ 2\ta-2(\ga-1) \}}\mathscr{E} (0).
\end{split}\ee
In a similar way to deriving \ef{sec41}, we have
\be\label{sec4.6}
x^3 \bar\rho^{2\ga-2\ta}(x)| r_x(x,t)-1|^2 \le C_{\iota}(1+t)^{-\frac{\beta}{2}+\frac{\beta}{2(\ga+\aa-\ta)}\max\{0, \ 4\ta-\ga-1\}}\mathscr{E} (0).
\ee
So,  \ef{deforvx} is a conclusion of  \ef{sec4.4}-\ef{sec4.6} and \ef{sec43}.

$$$$
\centerline{Acknowledgement}
Luo's research was  supported in part by NSF under grant DMS-1408839. Xin's research was partially supported by  the Zheng Ge Ru Foundation, and Hong
Kong RGC Earmarked Research Grants CUHK-4041/11P, CUHK-4048/13P, a Focus Area Grant from
The Chinese University of Hong Kong, and a grant from Croucher Foundation. Zeng's research  was supported in part by NSFC under grant 11301293, and the Center of Mathematical Sciences and Applications, Harvard University.

\newpage
\noindent {Tao Luo}\\
Department of Mathematics and Statistics\\
Georgetown University,\\
Washington, DC, 20057, USA. \\
Email: tl48@georgetown.edu

\vskip 0.25cm
\noindent Zhouping Xin\\
Institute of Mathematical Sciences\\
The Chinese University of Hong Kong\\
Shatin, NT, Hong Kong.\\
Email: zpxin@ims.cuhk.edu.hk

\vskip 0.25cm
\noindent Huihui Zeng\\
Yau Mathematical Sciences Center\\
Tsinghua University\\
Beijing, 100084, China;\\
Center of Mathematical Sciences and Applications\\
 Harvard University\\
  Cambridge, MA 02318, USA.\\
E-mail: hhzeng@mail.tsinghua.edu.cn

\end{document}